\documentclass[a4paper,12pt]{amsart}
\usepackage{amssymb}
\usepackage{ifthen}
\usepackage[dvips]{graphicx}
\usepackage{cite}
\nonstopmode \numberwithin{equation}{section}
\usepackage{amssymb}
\usepackage{ifthen}
\usepackage{graphicx}
\usepackage{amsmath}
\usepackage[T1]{fontenc} 
\usepackage[utf8]{inputenc}
\usepackage[usenames,dvipsnames]{color}
\usepackage{color}
\usepackage[english]{babel}
\usepackage{fancyhdr}
\usepackage{fancybox}
\usepackage{tikz}

\nonstopmode \numberwithin{equation}{section}
\theoremstyle{plain}

\newtheorem{prop}{Proposition}

\newtheorem{conj}{Conjecture}

\theoremstyle{definition}
\newtheorem{defn}{Definition}[section]
\newtheorem{example}{Example}[section]
\newtheorem{thm}{Theorem}[section]
\newtheorem{cor}{Corollary}[section]
\newtheorem{lem}{Lemma}[section]
\newtheorem{prob}{Problem}
\newtheorem{rem}{Remark}[section]
\newtheorem{ques}{Question}[section]

\newcounter{minutes}
\divide\time by 60
\newcounter{hours}
\multiply\time by 60
\addtocounter{minutes}{-\time}

\newcounter {own}
\def\theown {\thesection       .\arabic{own}}

\newenvironment{pf}[1][]{%
 \vskip 3mm
 \noindent
 \ifthenelse{\equal{#1}{}}%
  {{\slshape Proof. }}%
  {{\slshape #1.} }%
 }%
{\qed\bigskip}

\newcounter{alphabet}

\def\be{\begin{equation}}
\def\ee{\end{equation}}

\newcommand{\bee}{\begin{enumerate}}
\newcommand{\eee}{\end{enumerate}}

\newcommand{\blem}{\begin{lem}}
\newcommand{\elem}{\end{lem}}
\newcommand{\bthm}{\begin{thm}}
\newcommand{\ethm}{\end{thm}}
\newcommand{\bcor}{\begin{cor}}
\newcommand{\ecor}{\end{cor}}
\newcommand{\beg}{\begin{examp}}
\newcommand{\eeg}{\end{examp}}
\newcommand{\begs}{\begin{examples}}
\newcommand{\eegs}{\end{examples}}

\newcommand{\bdefn}{\begin{defn}}
\newcommand{\edefn}{\end{defn}}

\newcommand{\bprob}{\begin{prob}}
\newcommand{\eprob}{\end{prob}}
\newcommand{\bei}{\begin{itemize}}
\newcommand{\eei}{\end{itemize}}

\newcommand{\bcon}{\begin{conj}}
\newcommand{\econ}{\end{conj}}
\newcommand{\bcons}{\begin{conjs}}
\newcommand{\econs}{\end{conjs}}
\newcommand{\bprop}{\begin{prop}}
\newcommand{\eprop}{\end{prop}}
\newcommand{\br}{\begin{rem}}
\newcommand{\er}{\end{rem}}
\newcommand{\brs}{\begin{rems}}
\newcommand{\ers}{\end{rems}}
\newcommand{\bo}{\begin{obser}}
\newcommand{\eo}{\end{obser}}
\newcommand{\bos}{\begin{obsers}}
\newcommand{\eos}{\end{obsers}}
\newcommand{\bpf}{\begin{pf}}
\newcommand{\epf}{\end{pf}}
\newcommand{\ba}{\begin{array}}
\newcommand{\ea}{\end{array}}
\newcommand{\beq}{\begin{eqnarray}}
\newcommand{\beqq}{\begin{eqnarray*}}
\newcommand{\eeq}{\end{eqnarray}}
\newcommand{\eeqq}{\end{eqnarray*}}

\begin{document}

\title{Characterization of solutions of refined Fermat-type functional equations in $ \mathbb{C}^n $}

\author{Molla Basir Ahamed}
\address{Molla Basir Ahamed,
	Department of Mathematics,
	Jadavpur University,
	Kolkata-700032, West Bengal, India.}
\email{mbahamed.math@jadavpuruniversitry.in}

\author{Sanju Mandal}
\address{Sanju Mandal,
	Department of Mathematics,
	Jadavpur University,
	Kolkata-700032, West Bengal, India.}
\email{sanjum.math.rs@jadavpuruniversity.in}

\subjclass[{AMS} Subject Classification:]{Primary 39A45, 30D35, 35M30, 32W50}

\def\thefootnote{}
\footnotetext{ {\tiny File:~\jobname.tex,
printed: \number\year-\number\month-\number\day,
          \thehours.\ifnum\theminutes<10{0}\fi\theminutes }
} \makeatletter\def\thefootnote{\@arabic\c@footnote}\makeatother

\begin{abstract}
The main purpose of this article is concerned with the existence and the precise forms of the transcendental  solutions of several refined versions of Fermat-type functional equations with polynomial coefficients in several complex variables by utilizing the Nevanlinna theory of meromorphic functions in several complex variables. In fact, we investigate the existence and forms of the transcendental solutions of non-linear quadratic trinomial equations in $\mathbb{C}^n$. As a consequence of our result, we show that solutions of binomial equations in $\mathbb{C}^n$ can be explored and this exploration broaden the scope of the study of functional equations in $ \mathbb{C}^n $. The results we obtained are improvements over certain recent findings, noted as remarks. In addition, some examples relevant to the content of the paper have been exhibited in support of the validation of each results. Further, we discuss situations when solutions of such equations does not exist.
\end{abstract}

\maketitle
\pagestyle{myheadings}
\markboth{Molla Basir Ahamed  and  Sanju Mandal}{Characterization of solutions of Fermat-type functional equations in $ \mathbb{C}^n $}
\tableofcontents
\section{Introduction}
It is well known that, the complex oscillation theory of meromorphic solutions of differential equations is an important topic in complex analysis. Some results can be found in \cite{Laine-WG-1993}, where the Nevanlinna theory is an effective research tool. Recently, many results on meromorphic solutions of complex difference equations were rapidly obtained, such as \cite{Chiang & Feng & 2008, Chiang-Feng-TAMS-2009} and so on. In recent years, an extensive study of finding finite order entire and meromorphic solutions both in one and several complex variables of the functional equation 
\begin{align}\label{BS-eqq-1.1}
	f^n(z) + g^n(z)= 1\; \mbox{for integer}\; n\geq 1
\end{align}
what is called Fermat-type functional equations. In this study, Nevanlinna value distribution theory for meromorphic functions plays a significant role. In fact, the study of Fermat-type functional equations has a long history. In $ 1939 $, Iyer \cite{Iyer & J. Indian. Math. Soc. & 1939} studied the entire solutions of \eqref{BS-eqq-1.1} for $ n=2 $ and showed that the only solution is of the form $ f(z)=\cos(\varphi(z)) $ and $ g(z)=\sin(\varphi(z)) $, where $ \varphi $ is an entire function. Inspired by the results in \cite{Iyer & J. Indian. Math. Soc. & 1939}, many scholars paid considerable attention to investigating the existence and form of entire and meromorphic solutions of some variations of equations. We refer to the articles \cite{Xu-Liu-Li-JMAA-2020,Xu-RMJM-2021, XU-JIANG-RACSAM-2022} and references therein for some recent development on the topic.  In the past decade, with the establishment of difference analogues of Nevanlinna theory (see e.g., \cite{Cao-MN-2013, Chiang & Feng & 2008,Halburd & Korhonen & 2006}), many researchers studied the solutions of Fermat type difference and also Fermat-type differential-difference equations (see \cite{Liu-Song & Results Math & 2017,Qi-Liu-Yang-JCMA-2017,Haldar-MJM-2023}). Here, an equation is called differential-difference equation, if it includes derivatives, shifts or differences of $ f(z) $, which can be called $ DDE $ for short (see \cite{Lu-Lu-Li-Xu-RM-2019}). As far as we know, there are few results on solutions of complex partial differential-difference equations.\vspace{2mm}

The Fermat type equation \eqref{BS-eqq-1.1} appears in particle mechanics, in particular, it appears as the Lagrange function in the Lagrangian functional describing the “action” of a system. Finding the precise form of the solution with the help of Nevanlinna theory is an active area of research in function theory and researchers are paid considerable attention to this topic in recent years. For example, Liu \emph{et al.} \cite{Liu-Cao-Cao-AM-2012} have investigated the Fermat-type difference equation $ f^2 (z) + f^2 (z +c) = 1 $ in $ \mathbb{C} $ and obtained the transcendental entire solutions of finite order $ f(z)=\sin(Az+B) $, where $ B $ is a constant and $ A = ((4k+1)\pi)/(2c) $, where $ k $ is an integer. Later, Han and Lü \cite{Han-JCMA-Lu-2019} established the solution to the more general complex difference equation  $ f^n(z) +g^n(z) = e^{\alpha z+\beta} $. Moreover, Liu \emph{et al.} \cite{Liu-Cao-Cao-AM-2012} showed that the existence of solutions for the complex differential-difference equations $ {f^{\prime}}^2(z) + f(z +c)^2 = 1 $ and $ {f^{\prime}}^2(z) + (f(z +c) -f(z))^2 = 1 $ in $ \mathbb{C} $.\vspace{1.2mm}

We know that the partial differential equations (PDEs) are occurring in various areas of applied mathematics, such as fluid mechanics, nonlinear acoustics, gas dynamics, and traffic flow (see \cite{Courant-Hilbert-I-1962, Garabedian-W-1964}). In general, it is difficult to find entire and meromorphic solutions for a nonlinear PDE. By employing Nevanlinna theory and the method of complex analysis, there were a number of literature focusing on the solutions of some PDEs and their many variants, readers can refer to \cite{Cao-Xu-ADM-2020, Khavinson & Am. Math. Mon & 1995,Li-AM-2007, Gundersen-PEMS-2020, Saleeby-AM-2013}. The solutions of Fermat-type $ PDEs $ were investigated by \cite{Li-IJM-2004, Saleeby-A-1999,Haldar-MJM-2023}. Most noticeably, in $ 1995 $, Khavinson \cite{Khavinson & Am. Math. Mon & 1995} derived that any entire solution of the partial differential equation in $ \mathbb{C}^2 $,
\begin{align*}
	\left(\frac{\partial u}{\partial z_1}\right)^2+\left(\frac{\partial u}{\partial z_2}\right)^2=1
\end{align*}
is necessarily linear, i.e., $ u(z_1, z_2)=az_1+bz_2+c $, where $ a,b,c\in \mathbb{C} $, and $ a^2+b^2=1 $. This $ PDE $ in the real variable case occurs in the study of characteristic surfaces and wave propagation theory, and it is the two-dimensional eiconal equation, one of the main equations of geometric optics (see \cite{Courant-Hilbert-I-1962}). A more general equation of the form $F^2+G^2=e^{p(z)}$, where $p(z)$ is a polynomial has been studied recently to find the solutions. Due to the development of the logarithmic derivative lemma, in recent years, functional equations involving derivatives with shift or differences operators of meromorphic functions are studied extensively. Furthermore, Li \cite{Li-NMJ-2005, Li-AM-2007} have continued the research and discussed solutions of a series of $ PDEs $ with more general forms including $ \left({\partial f}/{\partial z_1}\right)^2 + \left({\partial f}/{\partial z_2}\right)^2 =e^g $, $ \left({\partial f}/{\partial z_1}\right)^2 + \left({\partial f}/ {\partial z_2}\right)^2 = p $, etc., where $ g, p $ are polynomials in $ \mathbb{C}^2 $. \vspace{1.2mm}

In \cite{Liu-Cao-Cao-AM-2012}, it is established that there is no finite order transcendental entire solution of Fermat-type difference equation $ f^n(z) + f^m(z+c) = 1 $ when $ n > m >1 $ or $ n= m >2 $. As a continuation of the research, later, the following result considering with polynomial coefficients in the equation in $ \mathbb{C} $ is established in \cite{Liu-Yang -CMFT-2013}.
\begin{thm}\label{th-1.1}\cite{Liu-Yang -CMFT-2013}
Let $P(z)$, $Q(z)$ be non-zero polynomials. If the difference equation
\begin{align}\label{E-1.2}
	f^2(z) + P^2(z)f^2(z+c)=Q(z)
\end{align}
admits a transcendental entire solution of finite order, then $P(z)\equiv \pm 1$ and $Q(z)$ reduces to a constant $q$. Thus, $f(z)=\sqrt{q} \sin(Az+B)$, where $B$ is a constant and $A={(4k+1)\pi}/{2c}$, where $k$ is a integer.
\end{thm}
In the next section, we give details on the study of fermat-type functional equations in several complex variables. 
\section{Description of the solutions of Fermat-type functional equations in $ \mathbb{C}^n $}
Over the past few years, Nevanlinna theory has made significant strides in the study of several complex variables. Notable advancements include the second fundamental theory of Nevanlinna and the application of difference analogues to logarithmic derivative lemmas, which have contributed substantially to the literature on holomorphic functions and their associated functional equations. Furthermore, through the application of the Nevanlinna theory for several complex variables with a difference perspective, Xu and Cao \cite{Xu-Cao-MJM-2018} investigated and characterized  the solution of the Fermat-type difference equation $ f^2(z) + f^2(z + c) = 1 $ in $ \mathbb{C}^n $. A significant characteristic of polynomials in $ \mathbb{C}^n $ (for $ n\geq 2 $) is their ability to exhibit periodicity, which is a fundamentally distinct behavior compared to polynomials in $ \mathbb{C} $. As a result, the solutions of functional equations in $ \mathbb{C}^n $ (for $ n\geq 2 $) with polynomial coefficients will possess distinct characteristics.\vspace{2mm} 

The study of solutions to Fermat-type equations in $\mathbb{C}^n$ is a captivating area of research that has recently garnered significant attention. Before delving into a detailed examination of such equations in $\mathbb{C}^n$, let us revisit a crucial result concerning Fermat-type equations in $\mathbb{C}^n$. In order to achieve symmetry, similar to the case of single-variable equations, we can summarize the results on Fermat-type functional equations in $\mathbb{C}^n$ for $n \geq 3$. By utilizing \cite[Theorems 1.1 and 1.2]{Saleeby-AM-2013}, along with the findings presented in \cite{Gross & Bull. Amer. & 1966, Gross-AMM-1966} and \cite{Baker-PAMS-1996}, we obtain the following theorem, which characterizes the entire and meromorphic solutions in $\mathbb{C}^n$ for the equation $f^m+g^m=1$ (possibly up to a fractional linear transformation of the uniformizing variable). It is important to mention 19th-century results on the global parameterization of elliptic curves, such as Picard's work (\cite{Saleeby-Analysis-1999}), Dixon's contributions to $sm$ and $cm$ functions (refer to \cite{Dixon-1989}), and Liouville's result stating that any elliptic function $f$ can be expressed as $F(\wp)+\wp^{\prime}G(\wp)$, where $\wp$ represents the Weierstrass elliptic function and $F$ and $G$ are rational functions of their arguments (see \cite{Neville-MG-1941}).

\begin{thm}\cite{Saleeby-AM-2013}
	For $ h:\mathbb{C}^n\rightarrow \mathbb{C} $ entire, the solutions of the equation $ f^m +g^m=1 $, $ m>1 $ are characterize as follows:
	\begin{enumerate}
		\item [(i)] for $ m=2 $, the entire solutions are $ f(z)=\cos(h(z)),\;g(z)=\sin(h(z))$;
		\item[(ii)] for $ m>2 $, there are no non-constant entire solutions;
		\item[(iii)] for $ m=2 $, the meromorphic solutions are $ f(z)={(1-\beta^2(z))}/{(1+\beta^2(z))} $ and $g(z)={2\beta(z)}/{(1+\beta^2(z))} $, with $ \beta $ being meromorphic on $ \mathbb{C}^n $.
		\item[(iv)] for $ m=3 $, the meromorphic solutions are of the form 
		\begin{align*}
			f(z)=\frac{1}{2\wp(h)}\left(1+\frac{\wp^{\prime}(h)}{\sqrt{3}}\right)\;\; \mbox{and}\;\; g(z)=\frac{\omega}{2\wp(h)}\left(1-\frac{\wp^{\prime}(h)}{\sqrt{3}}\right),
		\end{align*}
		where $ \left(\wp^{\prime}\right)^2=4\wp^3-1 $;
		\item[(v)] for $ m>3 $, there are no non-constant meromorphic solutions.
	\end{enumerate}
\end{thm}
With the introduction of the difference analogue of Nevanlinna theory for meromorphic functions with several complex variables, particularly the logarithmic derivative lemma proposed by Cao and Korhonen \cite{Cao-Korhonen-JMAA-2016}, numerous researchers have focused on exploring complex difference equations and complex Fermat-type difference equations with constant coefficients. This can be seen in various studies (e.g., \cite{Xu-RMJM-2021, Xu-AMP-2022} and their references). To gain a comprehensive understanding of Nevanlinna theory in $ \mathbb{C}^n $, we suggest referring to \cite{Ye-MZ-1996} and the literature mentioned therein. It is worth noting that non-constant polynomials in $ \mathbb{C}^n $ may exhibit periodicity when $ n\geq 2 $. Therefore, investigating solutions to Fermat-type difference equations in $ \mathbb{C}^n $ is an intriguing avenue of research, as these solutions possess fundamentally different characteristics compared to those in a single complex variable.\vspace{1.2mm}

Motivated by the above discussion and to continue the research for further investigations on Fermat-type functional equations in $\mathbb{C}^n$, we raised the following questions:
\begin{ques}\label{qu-1.1}
Can we derive a result corresponding to Theorem \ref{th-1.1} with polynomial coefficients in $ \mathbb{C}^n $?
\end{ques}
\begin{ques}\label{qu-1.2}
What can be said about the solutions of equations when $f(z)$ is replaced by $\frac{\partial f(z)}{\partial z_i}$ and $Q(z)$ is replaced by $Q(z)e^{g(z)}$ in Theorem \ref{th-1.1}, where $g$ is a polynomial in $ \mathbb{C}^n $?
\end{ques}
We now recall here an interesting result in $ \mathbb{C}^2 $, where an exploration presented in \cite{Zheng-Xu-AM-2022} delved into the existence and precise characteristics of transcendental entire solutions in $\mathbb{C}^2$.
\begin{thm}\cite{Zheng-Xu-AM-2022}\label{th-1.2}
Let $c=(c_1,c_2)\in\mathbb{C}^2\setminus\{(0,0)\}$ and $a_1, a_2, a_3$ be nonzero constant in $\mathbb{C}$. If the Fermat-type functional equation
\begin{align}\label{eq-1.2}
	a_1 f(z)^2 +(a_2 f(z+c) +a_3 f(z))^2 =1
\end{align}
has a transcendental entire solution $f(z)$ with finite order. Then $a^2_{1} +a^2_{3}=a^2_{2}$ and $f(z)$ is of the form
\begin{align*}
	f(z)=\dfrac{1}{a_1}\sin\left(L(z) +H(s) +B\right),
\end{align*}
where $L(z)=\alpha_1 z_1 +\alpha_2 z_2$, $\alpha_1,\alpha_2\in\mathbb{C}$,
$H(s)$ is a polynomial in $s:=c_2z_1 -c_1z_2$ in $\mathbb{C}^2$, and $L(z)$ satisfies $L(c)=\alpha_1 c_1 +\alpha_2 c_2=\theta+k\pi\pm {\pi}/{2}$, $\tan\theta={a_3}/{a_1}$.
\end{thm}
After in depth study of Theorem \ref{th-1.2} and its proof, we understand that two fold generalizations is possible. One is that the result can be generalized in $\mathbb{C}^n$ ($n\geq 2$) instead of $\mathbb{C}^2$, and the other one is that the functional equation in \eqref{eq-1.2} can be chosen with a more general setting. Henceforth, Theorem \ref{th-1.2} lead us to ask the following question.
\begin{ques}\label{Q-2.3}
What had happened about the solutions if the constant $1$ in R.H.S of \eqref{eq-1.2} is replaced by
$e^g$ in Theorem \ref{th-1.2}, where $g$ is a polynomial in $\mathbb{C}^n$ for $n\geq 2$?
\end{ques}
In this paper,  one of our aim is to answer Question \ref{Q-2.3}. Before delving into our investigations, we would like to highlight some results obtained by Saleeby regarding trinomial quadratic equations. In \cite{Saleeby-AM-2013}, Saleeby made important observations on right factors of meromorphic functions to describe complex analytic solutions to the quadratic functional equation $ f^2+2\alpha fg +g^2 =1 $, where $ \alpha^2\neq1 $ is a constant in $ \mathbb{C} $. Saleeby's work also associated with partial differential equations of the form $ u^2_x+2\alpha u_x u_y +u^2_y=1 $. Building upon this, Liu and Yang \cite{Liu & Yang & 2016} further investigated the existence and form of solutions for quadratic trinomial functional equations. It is proved in \cite{Liu & Yang & 2016} that if $ \alpha\neq \pm 1,0 $, then the equation $ f(z)^2+2\alpha f(z)f^{\prime}(z)+{f^{\prime}}^2(z)=1 $ has no transcendental meromorphic solutions. However, the equation $ f(z)^2+2\alpha f(z)f(z+c)+f(z+c)^2=1 $ must have transcendental entire solutions of order equal to one. In recent years, there has been extensive research on trinomial partial differential-difference equations for functions of two complex variables. Motivated by the techniques used by Saleeby in \cite{Saleeby-AM-2013} to solve trinomial equations, we raised the following questions corresponding to Theorems \ref{th-1.2}.
\begin{ques}\label{qu-1.3}
Do solutions exist for equations \eqref{eq-1.2} if the constant $1$ is replaced by a function $e^g$, where $g$ is a polynomial in $\mathbb{C}^n$? Moreover, with this setting, can we find solutions if the binomial equation is replaced by a quadratic trinomial equation in $\mathbb{C}^n$ with arbitrary coefficients?
\end{ques}\vspace{1.2mm}

Motivated by the above questions, this article is concerned with the description of entire solutions for several difference equations and partial differential-difference of Fermat type with more general forms in $\mathbb{C}^n$. Our main purpose of this paper is to discuss the finite order transcendental entire solutions of a Fermat-type binomial and quadratic trinomial functional equation with polynomial coefficients and find the precise form of the solutions. Henceforth, we consider the following equations
\begin{align}\label{Meq-1.1}
	af^2(z)+bP^2(z)[a_1f(z+c)+a_0f(z)]^2=Q(z)e^{g(z)},
\end{align}

\begin{align}\label{Meq-1.2}
	af^2(z+c)+bP^2(z)\left(\frac{\partial f(z)}{\partial z_i}\right)^2=Q(z)e^{g(z)}
\end{align}
and
\begin{align}\label{Meq-1.3}
	af^2(z) +2\omega f(z)\left(\gamma_1f(z+c)+\gamma_2 f(z)\right) +b\left(\gamma_1f(z+c)+\gamma_2 f(z)\right)^2=e^{g(z)},
\end{align}
where $ P(z) $, $ Q(z) $ and $ g(z) $ are non-zero polynomials in $ \mathbb{C}^n $. Moreover, we show that the precise form of the polynomial $ g(z) $ in $ \mathbb{C}^n $ can be obtained.\vspace{2mm}

This article is organized as follows:  In section 3, we  presents the main results of this paper, which address Fermat-type binomial and trinomial differential-difference equations with polynomial coefficients in $\mathbb{C}^n$. In fact, these main results correspond to the Questions \ref{qu-1.1}, \ref{qu-1.2}, \ref{Q-2.3}, and \ref{qu-1.3}, and we provide illustrative examples for each result to examine the precise form of transcendental entire solutions with finite order. We also include several remarks in support of our improved results.  In section 4, we provide a comprehensive review of essential Lemmas and present a detailed proof for each significant result.\vspace{2mm}

\section{Main results, corollaries and remarks}
The main tools used in this paper are the Nevanlinna theory and difference analogue of Nevanlinna theory with several complex variables. Throughout this paper, we assume that $ z+c=(z_1 + c_1, \ldots, z_n + c_n) $, for any $ z=(z_1,\ldots,z_n) $ and $ c=(c_1,\ldots,c_n) $ are in $ \mathbb{C}^n $. In fact, we know that any polynomial $ P(z) $ in $ \mathbb{C}^n $ of degree $ p $ can be expressed as $ P(z)=\sum_{|I|=0}^{p} a_{\alpha_1,\ldots,\alpha_n} z^{\alpha_1}_1\cdots z^{\alpha_n}_n $, where $ I=(\alpha_1,\ldots,\alpha_n) $  be a multi-index with $ |I|= \sum_{j=0}^{n}\alpha_j $ and $ \alpha_j $ are non-negative integers. For solving quadratic trinomial equations, we define 
\begin{align*}
\omega_1:=-\frac{\omega}{\sqrt{ab}} \pm \frac{\sqrt{\omega^2 -ab}}{\sqrt{ab}}\;\; \mbox{and}\;\;\omega_2: =-\frac{\omega}{\sqrt{ab}} \mp \frac{\sqrt{\omega^2 -ab}}{\sqrt{ab}}.
\end{align*}
Our key contribution in this paper pertains to Question \ref{qu-1.1} and demonstrates that the polynomial $P(z)$ discussed in \eqref{Meq-1.1} eventually attains a constant value, thereby enabling the determination of the solutions' precise form.
\begin{thm}\label{th-2.1}
Let $ c=(c_1,\ldots,c_n)\in\mathbb{C}^n\setminus\{(0,\ldots,0)\} $ and $ ab\neq 0 $. If $ f $ is a transcendental entire solution with finite order of the equation \eqref{Meq-1.1}, then the polynomial $ P(z) $ reduces to non-zero constant $ {\sqrt{a}}/{i\sqrt{b}\left(a_0 \pm a_1 e^{\frac{L_1(c)+ L_2(c)} {2}}\right)} $. Moreover, $ f $ and $ g $ must be one of the following to forms:
\vspace{1.2mm}
\begin{enumerate}
	\item [(i)] $ g(z)= \psi(s) + K+\sum_{i=1}^{n}(\alpha_i +\beta_i)z_i $, where $ \psi $ is a polynomial in $ s=\sum_{i=1}^{n}d_i z_i $ with $ \sum_{i=1}^{n}d_i c_i=0 $ with $ \psi(z+c) = \psi(z) $; Moreover,
	\begin{align*}
		f(z)=\dfrac{Q_1(z)e^{L_1(z) +\psi_1(s) + k_1}+Q_2(z)e^{L_2(z) +\psi_2(s) + k_2}}{2\sqrt{a}},
	\end{align*}
	where $ Q_1(z)= \sum_{i=1}^{n}m_{1i} z_i +H_1(s_1) +B_1 $ with $ \sum_{i=1}^{n}m_{1i} c_i= 0 $ and $ Q_2(z)= \sum_{i=1}^{n}m_{2i} z_i +H_2(s_1) +B_2 $ with $ \sum_{i=1}^{n}m_{2i} c_i= 0 $; $H_1(s_1), H_2(s_1)$ are polynomials in $s_1=\sum_{i=1}^{n}e_i z_i  $ with $\sum_{i=1}^{n}e_i c_i=0$,  $H_1(z+c)=H_1(z)$, $H_2(z+c)=H_2(z)$ and $B_1, B_2$ are constants. Furthermore, $ L_1(z)=\sum_{i=1}^{n} \alpha_i z_i$, $ L_2(z)=\sum_{i=1}^{n}\beta_i z_i$, $ \psi_1(s) $ and $ \psi_2(s) $ are polynomials in $ s $ in $ \mathbb{C}^n $ with $ \psi_1(z+c)= \psi_1(z) $ and $ \psi_2(z+c)= \psi_2(z) $, $k_1,k_2$ are constants.
	\item [(ii)] \begin{align*}
		f(z)=\beta(z)e^{L_{21}(z)}\; \mbox{and}\;g(z)=2 L_{21}(z) + B,
	\end{align*}
	where $ \beta(z) $ satisfies $ a\beta^2(z)+bP^2(z)\left(a_1\beta(z+c)+a_0\beta(z)\right)^2=e^{B}Q(z) $ and $ L_{21}(z)=\sum_{i=1}^{n}\gamma_i z_i $ and $ \gamma_1,\ldots,\gamma_n, B $ are constants
	\end{enumerate}
\end{thm}
Theorem \ref{th-2.1} yields several noteworthy observations, all of which we meticulously examine in the subsequent remark.
\begin{rem}
We note that 
\begin{enumerate}
	\item[(i)] Theorem \ref{th-2.1} improves that in \cite[Theorem 2.1, Theorem 2.3]{Liu-Cao-Cao-AM-2012} and \cite[Theorem 1.1] {Haldar-MJM-2023}, and \cite[Theorem 1.3, Theorem 1.5]{Chen-Hu-Zhang-JKMS-2016} in the sense that in this result equations are considered with polynomial coefficients in $ \mathbb{C}^n $.\vspace{1.5mm}
	
	\item [(ii)] In particular, if $n=1,m$; $ a= 1, b=1, a_1=1, a_0= 0 $ and $ g(z)= 2k\pi i $, $k\in\mathbb{Z}$, one can obtain the conclusions of Theorem \ref{th-1.1}\cite{Liu-Yang -CMFT-2013} in $\mathbb{C}$ and \cite[Theorem 1.1]{Haldar-MJM-2023} in $\mathbb{C}^n$. As a result, we obtain a more improved result.\vspace{1.5mm}
		
	\item[(iii)] In particular, if $ a_1= 1, a_0= -1 $ and $ g(z)= 2k\pi i $, $k\in\mathbb{Z}$, then it is easy to see that \eqref{Meq-1.1} reduces to the equation $ af^2(z)+bP^2(z)\left(\Delta_c f(z)\right)^2=Q(z) $, where $ P(z), Q(z) $ are two non-zero polynomials in $ \mathbb{C}^n $. Solving this equation, finally it can be  obtained that $bP^2(z)\left(Q(z+c)-Q(z)\right)=aQ(z)$. Clearly, if $\deg P(z)=p$ and $\deg Q(z)=q$, then we see that $2p+q-1=q$, we obtain $p={1}/{2}$, which is a contradiction. Therefore, this equation has no transcendental entire solutions with finite order. Consequently, our result is an improvement of that in \cite[Proposition 5.3]{Liu-JMAA-2009} and \cite[Theorem 2.3] {Liu-Cao-Cao-AM-2012}.\vspace{1.5mm}
	
	\item[(iv)] The conclusion of Theorem \ref{th-2.1} reveals also one crucial fact: if one considers the polynomial $ P(z) $ in equation \eqref{Meq-1.1} as non-constant, then \eqref{Meq-1.1} does not admit any solutions. 
\end{enumerate}
\end{rem}
\begin{rem}
In $\mathbb{C}$, it can be observed that when a polynomial $p(z)$ and a non-zero constant $c$ satisfy the equation $p(z+c) - p(z) = \xi$, with $\xi$ as a constant, the polynomial $p(z)$ is must be linear and can be written in the form $p(z) = az + b$, where $b$ represents a constant belonging to $\mathbb{C}$. But, if $ p(z) $ is a polynomial in $\mathbb{C}^n$ such that $p(z+c)-p(z)=\eta$, where $\eta$ is a constant in $\mathbb{C}$, then $p(z)$ must be of the form $p(z)=L(z)+H(s)+B$, where $L(z)=\sum_{i=1}^{n}a_n z_n$ and $H(s)$ is a polynomial in $s=\sum_{i=1}^{n}d_i z_i$ with $\sum_{i=1}^{n} d_i c_i=0$. In this case, it is worth pointing out that $p(z)$ is not necessarily linear, i.e., degree of $p(z)$ grater than or equal to $1$. Consequently, taking into account this property of a polynomial, it becomes clear that to verify the accuracy of this statement, the following examples demonstrate instances where the transcendental solutions of certain difference equations in $\mathbb{C}^n$ may order $\rho(f)\geq 1 $.
\end{rem}
\begin{example}
Let $c=(2,-1,3)\in\mathbb{C}^3$, the transcendental entire solutions in $ \mathbb{C}^3 $ of the difference equation
\begin{align*}
	7f^2(z)-\frac{{7}}{\left(5\pm\sqrt{3}e^{(4+i)+\frac{\pi}{2}}\right)^2}(&\sqrt{3}f(z+c)+5f(z))^2\\&\quad=Q_1(z)Q_2(z)e^{4z_1 +(2\pi +i)z_2 +(\pi +i)z_3 +\psi(s) +\frac{\pi i}{3}}
\end{align*}
where $Q_1(z)=(4z_1 +2z_2 -2z_3) +(2z_1+z_2-z_3)^5 + {\pi i}/{12}$, $Q_2(z)=(5z_1 +z_2 -3z_3) +(2z_1+z_2-z_3)^8 + {\pi i}/{13}$ and $ \psi(s) $ is a polynomial in $ s= z_1 +(2+3i)z_2 +iz_3 $; must be of the form
\begin{align*}
	f(z)&=\dfrac{1}{2\sqrt{7}}\left(Q_1(z)e^{3z_1 +iz_2 +\pi z_3+(z_1 +(2+3i)z_2 +iz_3)^{10} + \frac{\sqrt{3}\pi i}{7}} \right.\\&\quad\left.+ Q_2(z)e^{z_1 +2\pi z_2 +iz_3 +(z_1 +(2+3i)z_2 +iz_3)^7 + \frac{\sqrt{7}\pi}{11}}\right).
\end{align*}
It is easy to see that $\rho(f)=\max\{10,7\}=10>1$.
\end{example}
\begin{example}
Let $c=(5,-2,3)\in\mathbb{C}^3$, the transcendental entire solutions in $ \mathbb{C}^3 $ of the difference equation
\begin{align*}
	13f^2(z)-\frac{{13}}{\left(\frac{\sqrt{7}}{2}\pm\sqrt{5}e^{\frac{17-5i}{2}}\right)^2}\bigg(&\sqrt{5}f(z+c)+\frac{\sqrt{7}}{2}f(z)\bigg)^2\\&\quad=Q_1(z)Q_2(z)e^{5z_1 +(1 +i)z_2 -(2 +i)z_3 +\psi(s) +\frac{\sqrt{17}\pi i}{15}},
\end{align*}
where $Q_1(z)=(4z_1 +z_2 -6z_3)+(z_1 +z_2-z_3)^7 +{\pi i}/{4}$, $Q_2(z)=(6z_1 +3z_2 -8z_3)+(z_1 +z_2-z_3)^6 +{\pi i}/{9}$ and $ \psi(s) $ is a polynomial in $ s= 3z_1 +6z_2 -z_3 $; must be of the form
\begin{align*}
	f(z)&=\frac{1}{2\sqrt{13}}\left(Q_1(z)e^{2z_1 +z_2 -iz_3 +(3z_1 +6z_2 -z_3)^{13} + \frac{\sqrt{5}\pi i}{2}} \right.\\&\left.\quad +Q_2(z)e^{3z_1 +i z_2 -2z_3 +(3z_1 +6z_2 -z_3)^5 + \frac{\sqrt{7}\pi i}{5}}\right).
\end{align*}
 Clearly, we have $\rho(f)=\max\{13,5\}=13>1$.
\end{example}
The subsequent corollary for functional equations in $ \mathbb{C}^2 $, derived from Theorem \ref{th-2.1}, is as follows.
\begin{cor}
Let $ c=(c_1,c_2)\in\mathbb{C}^2\setminus\{(0,0)\} $ and $ a,b\neq 0 $. If $ f $ is a transcendental entire solution with finite order of the equation \eqref{Meq-1.1}, then the polynomials $ P(z) $ reduce to non-zero constant $ {\sqrt{a}}/{i\sqrt{b}\left(a_0 \pm a_1 e^{\frac{L_1(c)+ L_2(c)} {2}}\right)} $. Moreover, $ f $ and $ g $ must be one of the following to forms:
\vspace{1.2mm}
\begin{enumerate}
	\item [(i)] $ g(z)= (\alpha_1 +\beta_1)z_1 + (\alpha_2 +\beta_2)z_2 +\psi(s) + K $, where $ \psi $ is a polynomial in $ s=c_2 z_1 -c_1 z_2 $. Further,
	\begin{align*}
		f(z)=\dfrac{Q_1(z)e^{L_1(z) +\psi_1(s) + k_1}+Q_2(z)e^{L_2(z) +\psi_2(s) + k_2}}{2\sqrt{a}},
	\end{align*}
	where $ Q_1(z)= \sum_{i=1}^{2}m_{1i} z_i +H_1(s_1) +B_1 $, with $ \sum_{i=1}^{2}m_{1i} c_i= 0 $ and $ Q_2(z)= \sum_{i=1}^{2}m_{2i} z_i +H_2(s_1) +B_2 $ with $ \sum_{i=1}^{2}m_{2i} c_i= 0 $; where $H_1(s_1), H_2(s_1)$ is a polynomial is $s_1=c_2 z_1 -c_1 z_2$, $B_1, B_2$ are constants. Furthermore, $ L_1(z)=\sum_{i=1}^{2} \alpha_i z_i$, $ L_2(z)=\sum_{i=1}^{2}\beta_i z_i$, $ \psi_1(s) $ and $ \psi_2(s) $ are polynomials in $ s $ in $ \mathbb{C}^2 $, $k_1,k_2$ are constants.
	\item [(ii)] \begin{align*}
		f(z)=\beta(z)e^{L_{21}(z)}\; \mbox{and}\;g(z)=2 L_{21}(z) + B,
	\end{align*}
	where $ \beta(z) $ satisfies $ a\beta^2(z)+bP^2(z)\left(a_1\beta(z+c)+a_0\beta(z)\right)^2=e^{B}Q(z) $ and $ L_{21}(z)=\gamma_1 z_1 +\gamma_2 z_2 $ and $ \gamma_1, \gamma_2, B $ are constants.
\end{enumerate}
\end{cor}
The second result we obtain in this paper is the following, which answers to the Question \ref{qu-1.2}. We show that if $ f $ is a transcendental entire solution with finite order of the equation \eqref{Meq-1.2}, then the polynomials $ P(z) $ and $ Q(z) $ reduce to non-zero constants.
\begin{thm}\label{th-2.2}
Let $ c=(c_1,\ldots,c_n)\in\mathbb{C}^n\setminus\{(0,\ldots,0)\} $ and $ ab\neq 0 $. If $ f $ is a transcendental entire solution with finite order of the equation \eqref{Meq-1.2}, then the polynomials $ P(z) $ and $ Q(z) $ reduce to non-zero constants. Moreover, $ f $ and $ g $ must be one of the following forms:
\begin{enumerate}
	\item [(i)] \begin{align*}
		f(z)=\dfrac{\alpha_1 e^{h_1(z-c)} + \alpha_2 e^{h_2(z-c)}}{2\sqrt{a}}\; \mbox{and}\; g(z)=(r_1+r_2)+ \sum_{i=1}^{n}(a_i+ b_i)z_i
	\end{align*}
	where $ h_1(z)= r_1+\sum_{i=1}^{n}a_i z_i$ and $ h_2(z)=r_2+ \sum_{i=1}^{n}b_i z_i$, and satisfy the relations $ e^{\sum_{i=1}^{n}(a_i -b_i) c_i }=-a_i/b_i $; where $ r_1, r_2, a_i, b_i\in \mathbb{C} $ for $ i=1, 2, \ldots, n $. \vspace{1.2mm}
	\item [(ii)] 
	\begin{align*}
		f(z)= \gamma(z-c)e^{L_1(z-c)+ H(s)}\; \mbox{and}\; g(z)=2L_1(z) + 2H(s) + r_5,
	\end{align*}
	where $ L(z)= \sum_{i=1}^{n} \xi_i z_i $ and $ H(s) $ is a polynomial in $ s:=\sum_{i=1}^{n} d_i z_i $ in $ \mathbb{C}^n $ such that $ \sum_{i=1}^{n} d_i c_i= 0 $ with $ H(z+c)= H(z) $, $ \xi_1, \ldots, \xi_n, r_5  $ are constants and $ \gamma(z) $ satisfies 
	\begin{align*}
		a\gamma^2(z) + bP^2(z) e^{-2L_1(c)} \left(\frac{\partial \gamma(z-c)}{\partial z_i} + (\xi_i +d_i) \gamma(z-c)\right)^2 = e^{r_5}Q(z).
		\end{align*}
	\end{enumerate}
\end{thm}
\begin{rem}
We conclude some remark corresponding to Theorem \ref{th-2.2}.
\begin{enumerate}
	\item [(i)] We note that Theorem \ref{th-2.2} improves the results of \cite[Theorem 1.3]{Liu-Cao-Cao-AM-2012} and \cite[Theorem 1.3] {Haldar-MJM-2023}. \vspace{1.5mm}
	
	\item [(ii)] If $n=1$, $P(z)=1$, $Q(z)=1$, $ a= 1, b=1 $ and $ g(z)= 2k\pi i $, $k\in\mathbb{Z}$, then we obtain the conclusion of \cite[Theorem 1.3]{Liu-Cao-Cao-AM-2012}. Thus, our result exhibits a higher level of improvement compared to previous results. \vspace{1.5mm}
	
	\item [(iii)] In particular, if $n=2$, $ a= 1, b=1 $ and $ g(z)= 2k\pi i $, $k\in\mathbb{Z}$, one can obtain the conclusions of \cite[Theorem 1.3]{Haldar-MJM-2023} in $\mathbb{C}^2$. Accordingly, our result can be regarded as a more refined and improved version in comparison to the existing ones.\vspace{1.5mm}
\end{enumerate}
\end{rem}
Our third result of this paper is the following and we obtained the precise form of the solutions to a quadratic trinomial difference
equations in $\mathbb{C}^n$, which answers to the Question \ref{qu-1.3}. In fact, we obtain the following result under the assumption $ \omega^2\neq 0, ab $.
\begin{thm}\label{th-2.3}
Let $ c=(c_1,\ldots,c_n)\in\mathbb{C}^n\setminus\{(0,\ldots,0)\} $, $ ab\neq 0 $ and $\gamma_1, \gamma_2$ are nonzero constants, and $ \omega^2\neq 0, ab $. If the difference equation \eqref{Meq-1.3} admits a transcendental entire solution of finite order, then $ g(z) $ must be a polynomial of the form $ g(z)= L(z) + H(s) + B_3 $, where $ L(z) $ is a linear function of the form $ L(z)= a_1 z_1 +\cdots + a_n z_n $ and $ H(s) $ is a polynomial in $ s:= d_1 z_1 +\cdots + d_n z_n $ in $ \mathbb{C}^n $ with $ d_1 c_1 +\cdots + d_n c_n = 0 $ with $ H(z+c)= H(z) $, $ a_1,\ldots,a_n\in \mathbb{C} $. Further, $ f(z) $ must satisfy one of the following cases:
\begin{enumerate}
	\item [(i)] 
	\begin{align*}
		f(z) = \frac{\omega_2\xi^2 -\omega_1}{\xi\sqrt{a}(\omega_2 -\omega_1)} e^{\frac{1}{2}(L(z) + H(s) + B_3)},
	\end{align*}
	where
	\begin{align*}
		e^{\frac{a_1 c_1 +\cdots + a_n c_n}{2}} =\dfrac{(\omega_1 \gamma_2 \sqrt{b} -\sqrt{a}) -(\omega_2\gamma_2\sqrt{b} -\sqrt{a})\xi^2} {\gamma_1\sqrt{b}(\omega_2\xi^2 -\omega_1)}.
	\end{align*} 
	\item[(ii)] 
	\begin{align*}
		f(z) = \dfrac{\omega_2 e^{L_1(z) +H_1(s) + D_1} - \omega_1 e^{L_2(z) +H_2(s) + D_2}}{\sqrt{a}(\omega_2 -\omega_1)},
	\end{align*}
	where $ L_{l}(z)= a_{l1} z_1 +\cdots + a_{ln} z_n $ and $ H_{l}(s) $ for $ l=1, 2 $ are polynomial in $ s $ in $ \mathbb{C}^n $ and $ a_{l1},\ldots, a_{ln}, D_1, D_2\in\mathbb{C} $ for $ l=1, 2 $. Also, satisfying
	\begin{align*}
		L_1(z) +H_1(s)\neq L_2(z) +H_2(s), \;\; g(z) = L(z) + H(s) + D,
	\end{align*}
	where $ L(z) = L_1(z) + L_2(z) $, $ H(s) = H_1(s) + H_2(s) $ and $ D= D_1 + D_2 $ and 
	\begin{align*}
		-\frac{(\omega_2\gamma_2\sqrt{b}-\sqrt{a})}{\omega_2\gamma_1 \sqrt{b}} e^{-L_{1}(c)}\equiv 1 \;\;\;\;\mbox{and} \;\;\;\; -\frac{(\omega_1\gamma_2\sqrt{b}-\sqrt{a})}{\omega_1\gamma_1 \sqrt{b}} e^{-L_{2}(c)}\equiv 1.
	\end{align*}
\end{enumerate}
\end{thm}
\begin{rem}
Theorem \ref{th-2.3} strengthened the result for the quadratic trinomial equation with arbitrary coefficients of \cite[Theorem 2.2]{Zheng-Xu-AM-2022} in $ \mathbb{C}^n $.
\end{rem}
By the following examples, we show that forms of the solutions in part (i) and (ii) of Theorem \ref{th-2.3} are precise.
\begin{example}
For $ c=(5, 2, -3)\in\mathbb{C}^3 $, by a routine computation, it can be easily shown that the transcendental entire solutions in $ \mathbb{C}^3 $ of the difference equation	
\begin{align*}
	2f^2(z) +8 f(z)[5f(z+c)- 3 f(z)] &+ 3[5f(z+c)- 3 f(z)]^2\\&\quad =e^{3z_1 +\ln\left(\frac{10\pm 9\sqrt{10}}{5(4\pm 3\sqrt{10})}\right)z_2 +5z_3 +H(s) +\frac{5\pi i}{6}},
\end{align*}
where $ H(s) $ is a polynomial in $ s:= 4z_1 -z_2 +6z_3 $ of degree $ n\geq 1 $; must be of the form
\begin{align*}
	f(z_1,z_2,z_3) = \frac{(-4\mp 3\sqrt{10})}{\mp 4\sqrt{10}} e^{\frac{1}{2}\left(3z_1 +\ln\left(\frac{10\pm 9\sqrt{10}}{5(4\pm 3\sqrt{10})}\right)z_2 +5z_3 +(4z_1 -z_2 +6z_3)^n +\frac{5\pi i}{6}\right)},
\end{align*}
and the order of this transcendental entire solution is $\rho(f)=n\geq 1$.
\end{example}
\begin{example}
For $ c=( 2, -3, 1)\in\mathbb{C}^3 $, the transcendental entire solutions in $ \mathbb{C}^3 $ of the difference equation	
\begin{align*}
	3f^2(z) &+10 f(z)[3f(z+c)- 2 f(z)] +[3f(z+c)- 2 f(z)]^2\\&\quad =e^{11z_1 +\left(4-\frac{1}{3}\ln\left(\frac{7\mp 2\sqrt{22}}{3(5\mp \sqrt{22})}\right)\right)z_2 +\left(\ln\left(\frac{7\pm 2\sqrt{22}}{3(5\pm \sqrt{22})}\right)-10\right)z_3 +H(s)+\frac{181\pi i}{88}},
\end{align*}
where $ H(s) $ is a polynomial in $ s:= 5z_1 +4z_2 +2z_3 $; must be of the form
\begin{align*}
	f(z_1,z_2,z_3) &= \dfrac{(-5\mp\sqrt{22})}{\mp 2\sqrt{66}} e^{6z_1 +4z_2 +\ln\left(\frac{7\pm2\sqrt{22}}{3(5\pm\sqrt{22})}\right)z_3 +(5z_1 +4z_2 +2z_3)^n +\frac{7\pi i}{8}} \\&\quad - \dfrac{(-5\pm\sqrt{22})}{\mp 2\sqrt{66}} e^{5z_1 -\frac{1}{3}\ln\left(\frac{7\mp2\sqrt{22}}{3(5\mp\sqrt{22})}\right)z_2 -10z_3 + (5z_1 +4z_2 +2z_3)^m +\frac{13\pi i}{11}},
\end{align*}
and the order of this transcendental entire solution is $\rho(f)=\max\{n,m\}\geq 1$.
\end{example}
In Theorem \ref{th-2.3}, the complete solutions of \eqref{Meq-1.3} are thoroughly explored in the case when $ \omega^2\neq 0, ab $. However, the study of finding the solutions would be considered complete and comprehensive if we also determine the solutions of \eqref{Meq-1.3} when $ \omega^2= 0, ab $. Therefore, it is natural to answer the following question.
\begin{ques}\label{q-3.1}
	What can we say about the solutions of \eqref{Meq-1.3} when $ \omega=0 $ or $ \omega^2=ab $?
\end{ques}
We answer this Question \ref{q-3.1} in the subsequent discussion. However, as a consequence of Theorem \ref{th-2.3}, we obtain the result which encompasses the scope of Theorem \ref{th-1.2} in case of trinomial equations in $\mathbb{C}^2$.
\begin{cor}
Let $ c=(c_1,c_2)\in\mathbb{C}^2\setminus\{(0,0)\} $, and $\gamma_1, \gamma_2$ are nonzero constants and $ ab\neq 0 $, and $ \omega^2\neq ab $. If the difference equation
\begin{align*}
	af^2(z) +2\omega f(z)\left(\gamma_1f(z+c)+\gamma_2 f(z)\right) +b\left(\gamma_1f(z+c)+\gamma_2 f(z)\right)^2=e^{g(z)}
\end{align*}
admits a transcendental entire solution of finite order, then $ g(z) $ must be a polynomial function of the form $ g(z)= L(z) + H(s) + B_3 $, where $ L(z) $ is a linear function of the form $ L(z)= a_1 z_1 + a_2 z_2 $ and $ H(s) $ is a polynomial in $ s:= c_2 z_1-c_1 z_2 $ in $ \mathbb{C}^2 $. Further, $ f(z) $ must satisfy one of the following cases:
\begin{enumerate}
	\item [(i)] 
	\begin{align*}
		f(z) = \dfrac{\omega_2\xi^2 -\omega_1}{\xi\sqrt{a}(\omega_2 -\omega_1)} e^{\frac{1}{2}(L(z) + H(s) + B_3)},
	\end{align*}
	and satisfying
	\begin{align*}
		e^{\frac{a_1 c_1 + a_2 c_2}{2}} =\dfrac{(\omega_1 \gamma_2 \sqrt{b} -\sqrt{a}) -(\omega_2\gamma_2\sqrt{b} -\sqrt{a})\xi^2} {\gamma_1\sqrt{b}(\omega_2\xi^2 -\omega_1)}.
	\end{align*} 
	\item[(ii)] 
	\begin{align*}
		f(z) = \dfrac{\omega_2 e^{L_1(z) +H_1(s) + D_1} - \omega_1 e^{L_2(z) +H_2(s) + D_2}}{\sqrt{a}(\omega_2 -\omega_1)},
	\end{align*}
	where $ L_{l}(z)= a_{l1} z_1 + a_{l2} z_2 $ and $ H_{l}(s) $ for $ l=1, 2 $ are polynomial in $ s:= c_2 z_1-c_1 z_2 $ in $ \mathbb{C}^2 $ and $ D_1, D_2\in\mathbb{C} $. Also, satisfying
	\begin{align*}
		L_1(z) +H_1(s)\neq L_2(z) +H_2(s), \;\; g(z) = L(z) + H(s) + D,
	\end{align*}
	where $ L(z) = L_1(z) + L_2(z) $, $ H(s) = H_1(s) + H_2(s) $ and $ D= D_1 + D_2 $, and
	\begin{align*}
		-\dfrac{(\omega_2\gamma_2\sqrt{b}-\sqrt{a})}{\omega_2\gamma_1 \sqrt{b}} e^{-L_{1}(c)}\equiv 1 \;\;\;\;\mbox{and} \;\;\;\; -\dfrac{(\omega_1\gamma_2\sqrt{b}-\sqrt{a})}{\omega_1\gamma_1\sqrt{b}} e^{-L_{2}(c)}\equiv 1.
	\end{align*}
\end{enumerate}
\end{cor}
What could be the precise form of the solutions of \eqref{Meq-1.3} when $ \omega^2=ab $ or $ \omega^2=0 $ \emph{i.e., if} $ \omega=\pm \sqrt{a}\sqrt{b} $ or $ \omega=0 $, we now discuss in details in the following remark.
\begin{rem}\label{rem-3.5}
 For $ \omega=\pm \sqrt{a}\sqrt{b} $ and $ \omega=0 $, equation \eqref{Meq-1.3} can be respectively written as
\begin{align}\label{E-2.2}
	\sqrt{a}f(z)\pm \sqrt{b}\left(\gamma_1f(z+c)+\gamma_2 f(z)\right)=\pm e^{g(z)/2}
\end{align}
and 
\begin{align}\label{E-2.3}
	af^2(z) +b\left(\gamma_1f(z+c)+\gamma_2 f(z)\right)^2=e^{g(z)}.
\end{align}
It suffices to discuss the following three cases:\vspace{1.2mm}

\noindent{\bf Case 1}: Assume that $g(z)$ is a constant. However, in this case, it is easy to see that $e^{g(z)/2}$ is also a constant, say $E(\neq 0)$. Then by the standard method of finding solutions of complex difference equations with constant coefficients (see \cite{Ahamed-IJPAM-2022,Ahamed-RM-2019}), it can be shown that the precise form of the transcendental entire solutions in $\mathbb{C}^n$ of \eqref{E-2.2} will be of the form
\begin{align*}
	f(z)=\left(\mp\frac{(\sqrt{a}\pm\gamma_2\sqrt{b})}{\gamma_1\sqrt{b}}\right)^{z/c}\pi(z)\pm \frac{E}{\left(\sqrt{a}\pm(\gamma_1+ \gamma_2)\sqrt{b}\right)},
\end{align*}
where $ \pi $ is a transcendental entire function with $\pi(z+c)=\pi(z)$. \vspace{2mm}

\noindent{\bf Case 2}: If $g$ is non-constant and $ c $-periodic, \emph{i.e.} $ g(z+c)=g(z) $, then the precise form of the transcendental entire solutions in $\mathbb{C}^n$ of \eqref{E-2.2} will be of the form
\begin{align*}
	f(z)=\left(\mp\frac{(\sqrt{a}\pm\gamma_2\sqrt{b})}{\gamma_1\sqrt{b}}\right)^{z/c}\pi(z)\pm \frac{e^{g(z)/2}}{\left(\sqrt{a}\pm(\gamma_1+ \gamma_2)\sqrt{b}\right)},
\end{align*}
where $ \pi $ is a transcendental entire function with $ \pi(z+c)=\pi(z) $ and $ g $ is a polynomial with $ g(z+c)=g(z) $.\vspace{1.2mm}

\noindent{\bf Case 3}: If $g$ is non-constant and non-periodic. Then the precise form of the transcendental entire solutions in $\mathbb{C}^n$ of \eqref{E-2.2} will be of the form
\begin{align*}
	f(z)=\left(\mp\frac{(\sqrt{a}\pm\gamma_2\sqrt{b})}{\gamma_1\sqrt{b}}\right)^{z/c}\pi(z)\pm \frac{e^{g(z)/2}}{\left(\sqrt{a}\pm(\gamma_1+ \gamma_2)\sqrt{b}\right)},
\end{align*}
where $ \pi $ is a transcendental entire function with $ \pi(z+c)=\pi(z) $ and $ g(z) $ is a polynomial of the form $g(z)=L(z) +H_1(s)+B$, where $L(z)=\alpha_1 z_1 +\cdots+\alpha_n z_n$, $H_1(s_1)$ is a polynomial in $s_1=\sum_{i=1}^{n}e_i z_i$, with $\sum_{i=1}^{n}e_i c_i=0$, $H_1(z+c)=H_1(z)$ and $B$ is a constant, and satisfying
\begin{align*}
	\pm e^{\frac{L(c)}{2}}=\pm\left(\sqrt{a}\pm\sqrt{b(\gamma_1 +\gamma_2)}\right) \mp 1.
\end{align*}
\end{rem}
\subsection{Some examples in support of Remark \ref{rem-3.5}}
\begin{example}
	The transcendental entire solutions in $\mathbb{C}^n$ of the difference equation
	\begin{align*}
		3f(z)+ 2\left(\sqrt{3}f(z+c)+\sqrt{5} f(z)\right)= e^{g(z)/2}
	\end{align*}
	must be of the form
	\begin{align*}
		f(z)=\left(-\frac{(3+2\sqrt{5})}{2\sqrt{3}}\right)^{z/c}\sin\left(\dfrac{2\pi z}{c}\right)+ \frac{e^{g(z)/2}}{\left(3+2(\sqrt{5}+ \sqrt{3})\right)},
	\end{align*}
	where $ g $ is a polynomial with $ g(z+c)=g(z) $, \textit{i.e.,} $g$ is of the form $g(z)=\sum_{i=1}^{n}\gamma_i z_i +H(s)+A$, where $\sum_{i=1}^{n} \gamma_i c_i=0$, $H(s)$ is a polynomial is $s=\sum_{i=1}^{n}d_i z_i$ with $\sum_{i=1}^{n}d_i c_i=0$, $H(z+c)=H(z)$ and $A$ is a constant.
\end{example}

\begin{example}
	The transcendental entire solutions in $\mathbb{C}^n$ of the difference equation
	\begin{align*}
		\sqrt{7}f(z)+ \sqrt{3}\left(2f(z+c)-f(z)\right)= -e^{g(z)/2}
	\end{align*}
	must be of the form
	\begin{align*}
		f(z)=\left(-\frac{(\sqrt{7}-\sqrt{3})}{2\sqrt{3}}\right)^{z/c}\cos\left(\dfrac{2\pi z}{c}\right)- \frac{e^{g(z)/2}}{\left(\sqrt{7} +\sqrt{3} \right)},
	\end{align*}
	where $ g $ is a polynomial with $ g(z+c)=g(z) $, \textit{i.e.,} $g$ is of the form $g(z)=\sum_{i=1}^{n}\gamma_i z_i +H(s)+A$, where $\sum_{i=1}^{n} \gamma_i c_i=0$, $H(s)$ is a polynomial is $s=\sum_{i=1}^{n}d_i z_i$ with $\sum_{i=1}^{n}d_i c_i=0$, $H(z+c)=H(z)$ and $A$ is a constant.
\end{example}

\begin{example}
	The transcendental entire solutions in $\mathbb{C}^n$ of the difference equation
	\begin{align*}
		\sqrt{5}f(z)- \sqrt{2}\left(3f(z+c)-2f(z)\right)= e^{g(z)/2}
	\end{align*}
	must be of the form
	\begin{align*}
		f(z)=\left(\frac{(\sqrt{5}+2\sqrt{2})}{3\sqrt{2}}\right)^{z/c}\sin\left(\dfrac{2\pi z}{c}\right)+ \frac{e^{g(z)/2}}{\left(\sqrt{5} -\sqrt{2} \right)},
	\end{align*}
	where $ g $ is a polynomial with $ g(z+c)=g(z) $, \textit{i.e.,} $g$ is of the form $g(z)=\sum_{i=1}^{n}\gamma_i z_i +H(s)+A$, where $\sum_{i=1}^{n} \gamma_i c_i=0$, $H(s)$ is a polynomial is $s=\sum_{i=1}^{n}d_i z_i$ with $\sum_{i=1}^{n}d_i c_i=0$, $H(z+c)=H(z)$ and $A$ is a constant.
\end{example}

\begin{example}
	The transcendental entire solutions in $\mathbb{C}^n$ of the difference equation
	\begin{align*}
		\sqrt{2}f(z)- \sqrt{3}\left(\sqrt{7}f(z+c)+\sqrt{11} f(z)\right)= -e^{g(z)/2}
	\end{align*}
	must be of the form
	\begin{align*}
		f(z)=\left(\frac{(\sqrt{2}-\sqrt{33})}{\sqrt{21}}\right)^{z/c}\cos\left(\dfrac{2\pi z}{c}\right)- \frac{e^{g(z)/2}}{\left(\sqrt{2} -(\sqrt{7}+ \sqrt{11})\sqrt{3} \right)},
	\end{align*}
	where $ g $ is a polynomial with $ g(z+c)=g(z) $, \textit{i.e.,} $g$ is of the form $g(z)=\sum_{i=1}^{n}\gamma_i z_i +H(s)+A$, where $\sum_{i=1}^{n} \gamma_i c_i=0$, $H(s)$ is a polynomial is $s=\sum_{i=1}^{n}d_i z_i$ with $\sum_{i=1}^{n}d_i c_i=0$, $H(z+c)=H(z)$ and $A$ is a constant.
\end{example}
It is worth noticing that if we allow $ \omega=0 $ in \eqref{Meq-1.3}, then it reduces to binomial equation \eqref{E-2.3} of the form $ F^2+G^2=e^g $. We need to explore the solutions in this case too. Hence, we obtain the following result for the case $\omega=0$, \textit{i.e.,} corresponding to the equation \eqref{E-2.3}. Consequently, the following result is an answer to Question \ref{Q-2.3}
\begin{thm}\label{th-2.4}
Let $ c=(c_1,\ldots,c_n)\in\mathbb{C}^n\setminus\{(0,\ldots,0)\} $, $ ab\neq 0 $ and $\gamma_1, \gamma_2$ are nonzero constants. If the difference equation \eqref{E-2.3} admits a transcendental entire solution of finite order, then $ g(z) $ must be a polynomial of the form $ g(z)= L(z) + H(s) + A $, where $ L(z) $ is a linear function of the form $ L(z)= a_1 z_1 +\cdots + a_n z_n $ and $ H(s) $ is a polynomial in $ s:= d_1 z_1 +\cdots + d_n z_n $ in $ \mathbb{C}^n $ with $ d_1 c_1 +\cdots + d_n c_n = 0 $ with $ H(z+c)= H(z) $, $ a_1,\ldots,a_n, A\in \mathbb{C} $. Further, $ f(z) $ must satisfy one of the following cases:
\begin{enumerate}
	\item [(i)] 
	\begin{align*}
		f(z) = \dfrac{(\xi^2 + 1)}{2\xi\sqrt{a}} e^{\frac{1}{2}(L(z) + H(s) +A)},
	\end{align*}
	where $ g(z)= L(z) + H(s) +A $, where $ L(z)= a_1 z_1 +\cdots + a_n z_n $; $ a_1,\ldots,a_n, A\in\mathbb{C} $ satisfying
	\begin{align*}
		e^{\frac{L(c)}{2}} =\dfrac{\sqrt{a}(\xi^2 -1) -i\gamma_2\sqrt{b}(\xi^2 +1)} {i\gamma_1\sqrt{b}(\xi^2 +1)}.
	\end{align*} 
	\item[(ii)] 
	\begin{align*}
		f(z) = \dfrac{e^{L_1(z) +H_1(s) + B_1} + e^{L_2(z) +H_2(s) + B_2}}{2\sqrt{a}},
	\end{align*}
	where $ L_{l}(z)= a_{l1} z_1 +\cdots + a_{ln} z_n $ and $ H_{l}(s) $ for $ l=1, 2 $ are polynomial in $ s $ in $ \mathbb{C}^n $ and $ a_{l1},\ldots, a_{ln}, B_1, B_2\in\mathbb{C} $ for $ l=1, 2 $. Also, satisfying
	\begin{align*}
		L_1(z) +H_1(s)\neq L_2(z) +H_2(s), \;\; g(z) = L(z) + H(s) +B,
	\end{align*}
	where $ L(z) = L_1(z) + L_2(z) $, $ H(s) = H_1(s) + H_2(s) $ and $ B= B_1 + B_2 $ and 
	\begin{align*}
		\dfrac{(i\gamma_2\sqrt{b}-\sqrt{a})}{-i\gamma_1 \sqrt{b}} e^{-L_{1}(c)}\equiv 1 \;\;\;\;\mbox{and} \;\;\;\; \dfrac{(i\gamma_2\sqrt{b}+\sqrt{a})}{-i\gamma_1 \sqrt{b}} e^{-L_{2}(c)}\equiv 1.
	\end{align*}
\end{enumerate}
\end{thm}
An interesting observation is that, in particular, if $\sqrt{a}=a_1, \sqrt{b}\gamma_1=a_2,\sqrt{b}\gamma_2=a_3$ and $ g(z)= 2k\pi i $, $k\in\mathbb{Z}$ in \eqref{E-2.3}, one can obtain the conclusions of \cite[Theorem 2.2]{Zheng-Xu-AM-2022} in $\mathbb{C}^2$. Accordingly, our result can be regarded as a more refined and improved version in comparison to \cite[Theorem 2.2]{Zheng-Xu-AM-2022}. 
\begin{cor}
Let $ c=(c_1,c_2)\in\mathbb{C}^2\setminus\{(0,0)\} $, $ ab\neq 0 $ and $\gamma_1, \gamma_2$ are nonzero constants. If the difference equation \eqref{E-2.3} admits a transcendental entire solution of finite order, then $ g(z) $ must be a polynomial of the form $ g(z)= L(z) + H(s) + A $, where $ L(z) $ is a linear function of the form $ L(z)= a_1 z_1 + a_2 z_2 $ and $ H(s) $ is a polynomial in $ s:= c_2 z_1 - c_1 z_2 $ in $ \mathbb{C}^2 $ with $ H(z+c)= H(z) $, $ a_1,a_2, A\in \mathbb{C} $. Further, $ f(z) $ must satisfy one of the following cases:
\begin{enumerate}
	\item [(i)] 
	\begin{align*}
		f(z) = \dfrac{(\xi^2 + 1)}{2\xi\sqrt{a}} e^{\frac{1}{2}(L(z) + H(s) +A)},
	\end{align*}
	where
	\begin{align*}
		e^{\frac{L(c)}{2}} =\dfrac{\sqrt{a}(\xi^2 -1) -i\gamma_2\sqrt{b}(\xi^2 +1)} {i\gamma_1\sqrt{b}(\xi^2 +1)}.
	\end{align*} 
	\item[(ii)] 
	\begin{align*}
		f(z) = \dfrac{e^{L_1(z) +H_1(s) + B_1} + e^{L_2(z) +H_2(s) + B_2}}{2\sqrt{a}},
	\end{align*}
	where $ L_{l}(z)= a_{l1} z_1 + a_{l2} z_2 $ and $ H_{l}(s) $ for $ l=1, 2 $ are polynomial in $ s $ in $ \mathbb{C}^2 $ and $ a_{l1},\ldots, a_{ln}, B_1, B_2\in\mathbb{C} $ for $ l=1, 2 $, satisfying
	\begin{align*}
		L_1(z) +H_1(s)\neq L_2(z) +H_2(s), \;\; g(z) = L_3(z) + H_3(s) +B_3,
	\end{align*}
	where $ L_3(z) = L_1(z) + L_2(z) $, $ H_3(s) = H_1(s) + H_2(s) $ and $ B_3= B_1 + B_2 $ and 
	\begin{align*}
		\frac{(i\gamma_2\sqrt{b}-\sqrt{a})}{-i\gamma_1 \sqrt{b}} e^{-L_{1}(c)}\equiv 1 \;\;\;\;\mbox{and} \;\;\;\; \frac{(i\gamma_2\sqrt{b}+\sqrt{a})}{-i\gamma_1 \sqrt{b}} e^{-L_{2}(c)}\equiv 1.
	\end{align*}
\end{enumerate}
\end{cor}
\section{Key lemmas and Proof of the main results}
We shall use the following standard notations of value distribution theory (see details \cite{Hu-Li-Yang-2003,Yang-Yi-2006}): $ T(r,f),\; m(r,f),\; N(r,f),\; \overline{N}(r,f),\ldots $ in $ \mathbb{C}^n $. We denote by $ S(r, f) $, any function satisfying $ S(r, f) = \circ\{T(r,f)\} $  as $ r\rightarrow\infty $, possibly outside a set of finite measure. Define the order of $ f $ by
\begin{align*}
	\rho(f):=\overline{\lim_{r\rightarrow\infty}}\dfrac{\log^{+}T(r,f)}{\log r}.
\end{align*}
First, we recall here some necessary lemmas which will play key roles in proving the main results.
\begin{lem}\cite{Ronkin_AMS-1974,Stoll-AMS-1974}\label{lem-3.1}
	For any entire function $ F $ on $ \mathbb{C}^n $, $ F(0)\neq 0 $ and put $\rho(n_F)=\rho < \infty $, where $ \rho(n_F) $ denotes be the order of the counting function of zeros of $ F $. Then there exist a canonical function $ f_F $ and a function $ g_F \in\mathbb{C}^n $ such that $ F(z) =f_F (z)e^{g_F (z)} $. For the special case $ n = 1 $, $ f_F $ is the canonical product of Weierstrass.
\end{lem}
\begin{lem}\cite{Polya-JLMS-1926}\label{lem-3.2}
	If $  g $ and $ h $ are entire functions on the complex plane $ \mathbb{C} $ and $ g(h) $ is an entire function of finite order, then there are only two possible cases: either
	\begin{enumerate}
		\item [(i)] the internal function $ h $ is a polynomial and the external function $ g $ is of finite order; or
		\item[(ii)] the internal function $ h $ is not a polynomial but a function of finite order, and the external function $ g $ is of zero order.
	\end{enumerate}
\end{lem}
\begin{lem}\cite{Hu-Li-Yang-2003}\label{lem-3.3}
	Let $ f_j(\not \equiv 0) \;(j=1,2,\ldots, m;\;m\geq 3) $ be meromorphic functions on $ \mathbb{C}^n $ such that $ f_1,\ldots,f_{m-1} $ are non-constant and $ f_1 +\ldots + f_m = 1 $, and such that 
	\begin{align*}
		\bigg|\bigg|\;\sum_{k=1}^{m}\left\{N_{m-1}\left(r,\frac{1}{f_k}\right) + (m-1)\overline{N}(r,f_k)\right\}<\lambda T(r,f_j) + O(\log^{+} T(r,f_j)),
	\end{align*}
	holds for $ j=1,\ldots, m-1 $, where $ \lambda<1 $ is a positive number. Then  $ f_m = 1 $.
\end{lem}

\begin{lem}\cite{Hu-Li-Yang-2003}\label{lem-3.4}
Suppose that $ a_0(z), a_1(z),\ldots, a_m(z) \; (m\geq 1) $ are meromorphic functions on $ \mathbb{C}^n $ and $ g_0(z), g_1(z),\ldots, g_m(z) $ are entire functions on $ \mathbb{C}^n $ such that $ g_i(z)- g_j(z) $ are not constants for $ 0\leq i< j\leq m $. If
\begin{align*}
	\sum_{i=0}^{m} a_i(z) e^{g_i(z)} \equiv 0
\end{align*}
and $ || T(r, a_i) = o(T(r)), \; i=0, 1, \ldots, m $ holds, where $ T(r):=\displaystyle\min_{0\leq i< j\leq m} T(r, e^{g_i -g_j}) $, then $ a_i(z) \equiv 0 $ for $ i=0, 1, \ldots, m $.
\end{lem}
For the convenience of the reader, we will present our proof of Theorems \ref{th-2.1}, \ref{th-2.2}, and \ref{th-2.3} in all detail.
\subsection{Proof of Theorem \ref{th-2.1}}
\begin{proof}[\bf Proof of Theorem \ref{th-2.1}] For the sake of simplicity, first we define $ L_c(f(z)):=a_1 f(z+c)+a_0f(z) $, where $ a_1, a_0\in\mathbb{C} $. It is easy to see that the given equation \eqref{Meq-1.1} can be expressed as 
	\begin{align}\label{eq-3.1}
		\left(\frac{\sqrt{a}f(z)+i\sqrt{b}P(z)\left(L_c(f(z))\right)}{e^{\frac{g(z)}{2}}}\right)\left(\frac{\sqrt{a}f(z)-i\sqrt{b}P(z)\left(L_c(f(z))\right)}{e^{\frac{g(z)}{2}}}\right)=Q(z).
	\end{align}
	In view of \eqref{eq-3.1} and Lemma \ref{lem-3.1}, the following pair of equations can be obtained
	\begin{align}\label{eq-3.2}
		\begin{cases}
			\sqrt{a}f(z)+iP(z)\left(a_1f(z+c)+a_0f(z)\right)=Q_1(z)e^{\frac{g(z)}{2}+\alpha(z)} \vspace{1.2mm}\\\sqrt{a}f(z)-iP(z)\left(a_1f(z+c)+a_0f(z)\right) =Q_2(z)e^{\frac{g(z)}{2}-\alpha(z)},
		\end{cases}
	\end{align}
	where $ \alpha(z) $ is an non-constant entire function on $ \mathbb{C}^n $ and $Q(z)=Q_1(z)Q_2(z)$. For brevity, we define
	\begin{align}\label{eq-3.3}
		h_1(z):=\dfrac{g(z)}{2}+\alpha(z) \;\mbox{and}\; h_2(z):=\dfrac{g(z)}{2}-\alpha(z).
	\end{align} 
	Since $ f $ is a transcendental entire function with finite order, in view of Lemma
	\ref{lem-3.2}, we see that $ \alpha(z) $ must be a polynomial in $ \mathbb{C}^n $. Solving the pair of equations in \eqref{eq-3.2}, we easily obtain
	\begin{align}\label{eq-3.4}
		\begin{cases}
			f(z)=\dfrac{Q_1(z)e^{h_1(z)}+Q_2(z)e^{h_2(z)}}{2\sqrt{a}}\vspace{2mm}\\
			a_1f(z+c)+a_0f(z)=\dfrac{Q_1(z)e^{h_1(z)}-Q_2(z)e^{h_2(z)}}{2i\sqrt{b}P(z)}.
		\end{cases}
	\end{align} 
	Combining both the expressions in \eqref{eq-3.4}, a simple computation shows that
	\begin{align}\label{eq-3.5}
		\dfrac{a_1}{\sqrt{a}}\dfrac{Q_1(z+c)}{P_1(z)Q_1(z)}e^{h_1(z+c)-h_1(z)} + \dfrac{a_1}{\sqrt{a}}\dfrac{Q_2(z+c)}{P_1(z)Q_1(z)}e^{h_2(z+c)-h_1(z)}-\dfrac{Q_2(z)}{Q_1(z)}e^{h_2(z)-h_1(z)}=1,
	\end{align} 
	where $ P_1(z)=\frac{1}{i\sqrt{b}P(z)} - \frac{a_0}{\sqrt{a}} $.\vspace{1.2mm}
	
	Moreover, the equation \eqref{eq-3.5} can be written as $g_1(z)+ g_2(z) +g_3(z)= 1$, where
	\begin{align*}
		\begin{cases}
			g_1(z)= \dfrac{a_1}{\sqrt{a}}\dfrac{Q_1(z+c)}{P_1(z)Q_1(z)} e^{h_1(z+c)-h_1(z)} \vspace{2mm}\\ g_2(z)= \dfrac{a_1} {\sqrt{a}}\dfrac{Q_2(z+c)}{P_1(z)Q_1(z)}e^{h_2(z+c)-h_1(z)} \vspace{2mm}\\ g_3(z)= -\dfrac{Q_2(z)} {Q_1(z)} e^{h_2(z)-h_1(z)}
		\end{cases}
	\end{align*}
	 In order to complete the proof, it is enough to discuss the following two possible cases. \vspace{1.5mm}
	
	\noindent{\bf Case 1.} Suppose that $ e^{h_2(z)-h_1(z)} $ is not a constant. Then, it is easy to see that polynomial $ h_2(z)-h_1(z) $ is a non-constant. Our aim is to show that both the functions  $ g_2 $ and $ g_3 $ are non-constants. Otherwise, suppose that $ g_2(z)=D_1 $ and $ g_3(z)=-D_2 $, where $ D_1, D_2 $ are both constants. Then, an easy computation leads to
	\begin{align}\label{eq-3.6}
		e^{h_2(z+c)-h_1(z)}=\frac{D_1\sqrt{a}}{a_1}\frac{P_1(z)Q_1(z)}{Q_2(z+c)} \;\mbox{and}\; e^{h_2(z)-h_1(z)}=\frac{D_2 Q_1(z)}{Q_2(z)}.
	\end{align}
	It is easy to see from \eqref{eq-3.6} that the left side is a transcendental entire function, whereas the right side of it is a rational function, which is a contradiction. Consequently, we see that
	\begin{align*}
		\sum_{k=1}^{3}\left\{N_{2}\left(r,\frac{1}{g_k}\right) + 2\overline{N}(r,g_k)\right\}= O(\log^{+} T(r,g_3))
	\end{align*}
	as $r$ sufficiently large outside possibly a set $E$ of $r$ with finite logarithmic measure. Thus, in view of Lemma \ref{lem-3.3} (with $m=3$), we must have $g_1(z)\equiv 1$ \textit{i.e.,}
	\begin{align}\label{eq-3.7}
		a_1 Q_1(z+c)e^{h_1(z+c)-h_1(z)}\equiv\sqrt{a} P_1(z) Q_1(z).
	\end{align}
	We observe that left side of \eqref{eq-3.7} is transcendental entire, whereas right
	side is polynomial, and this shows that $ h_1(z+c)-h_1(z)=\xi $, where $ \xi\in\mathbb{C} $ is a constant. An easy computation shows that $ h_1(z)=L_1(z) +\psi_1(s) + k_1 $, where $ L_1(z)=\alpha_1 z_1 +\cdots +\alpha_n z_n $ and $ \psi_1(s) $ is a polynomial in $ s:= d_1 z_1 +\cdots + d_n z_n $ in $ \mathbb{C}^n $ with $ d_1 c_1 +\cdots + d_n c_n = 0 $ and $ \psi_1(z+c)= \psi_1(z) $; $ \alpha_1,\ldots,\alpha_n, k_1\in \mathbb{C} $.\vspace{1.2mm} 
	
	Further, from \eqref{eq-3.7}, we see that
	\begin{align}\label{eq-3.8}
		a_1 Q_1(z+c)e^{L_1(c)}\equiv\sqrt{a} P_1(z) Q_1(z).
	\end{align}
	Since $Q_1(z)$ is a non-zero polynomial, from \eqref{eq-3.8}, we must have
	\begin{align}\label{Eq-4.9}
		P_1(z)=\frac{a_1 e^{L_1(c)}}{\sqrt{a}}.
	\end{align} 
	Thus, it follows from \eqref{eq-3.8} that $ Q_1(z+c)=Q_1(z) $ \emph{i.e.} the polynomial $ Q_1 $  is $ c $-periodic. By a simple computation, it can be easily shown that $ Q_1(z)= m_{11} z_1+ \cdots+ m_{1n} z_n +H_1(s_1)+B_1 $, where $ m_{11} c_1+\cdots +m_{1n} c_n= 0 $, $H_1(s_1)$ is a polynomial is $s_1=\sum_{i=1}^{n}e_i z_i  $ with $\sum_{i=1}^{n}e_i c_i=0$,  $H_1(z+c)=H_1(z)$ and $B_1$ is a constant.\vspace{1.2mm}
	
	Therefore, in view of \eqref{eq-3.5} and \eqref{eq-3.7}, we obtain
	\begin{align}\label{eq-3.9}
		a_1 Q_2(z+c) e^{h_2(z+c)-h_2(z)}=\sqrt{a}P_1(z)Q_2(z).
	\end{align}
     By the similar argument being used for $h_1(z)$, we obtain $ h_2(z)=L_2(z) +\psi_2(s) + k_2 $, where $ L_2(z)=\beta_1 z_1 +\cdots+ \beta_n z_n $ and $ \psi_2(s) $ is a polynomial in $ s:= d_1 z_1 +\cdots + d_n z_n $ in $ \mathbb{C}^n $ with $ d_1 c_1 +\cdots + d_n c_n = 0 $ with $ \psi_2(z+c)= \psi_2(z) $; $ \beta_1,\ldots,\beta_n, k_2\in \mathbb{C} $. Hence, the equation \eqref{eq-3.9} reduces to
	\begin{align}\label{Eq-4.10}
		a_1 Q_2(z+c) e^{L_2(c)}=\sqrt{a}P_1(z)Q_2(z).
	\end{align}
	Since $Q_2(z)$ is a non-zero polynomial, \eqref{eq-3.9} shows that
	\begin{align}\label{Eq-4.12}
		P_1(z)=\frac{a_1 e^{L_2(c)}}{\sqrt{a}}.
	\end{align}
	 As a matter of fact, from \eqref{Eq-4.10}, we see that $ Q_2(z+c)=Q_2(z) $. This shows that $ Q_2(z)= m_{21} z_1+ \cdots+ m_{2n} z_n +H_2(s_1)+B_2 $, where $ m_{21} c_1+\cdots +m_{2n} c_n= 0 $, $H_2(s_1)$ is a polynomial is $s_1=\sum_{i=1}^{n}e_i z_i  $ with $\sum_{i=1}^{n}e_i c_i=0$, $H_2(z+c)=H_2(z)$ and $B_2$ is a constant. Multiplying \eqref{Eq-4.9} and \eqref{Eq-4.12}, it is easy to see that 
	\begin{align}\label{Eq-4.13}
		aP^2_{1}(z)= a^2_{1} e^{L_1(z)+L_2(c)}.
	\end{align}
	Now plugging $P_1(z)=\frac{1}{i\sqrt{b}P(z)} - \frac{a_0}{\sqrt{a}}$ into \eqref{Eq-4.13}, we see that $P(z)$ becomes a constant of the form
	\begin{align*}
		P(z) =\dfrac{\sqrt{a}}{i\sqrt{b}\left(a_0 \pm a_1 e^{\frac{L_1(c)+ L_2(c)} {2}}\right)}= \mbox{constant}.
	\end{align*}
	Moreover, from  \eqref{eq-3.3}, it is easy to see that
	\begin{align*}
		g(z)= L(z) + \psi(s) + K= (\alpha_1 +\beta_1)z_1 + \cdots + (\alpha_n +\beta_n)z_n + \psi(s) + K,
	\end{align*}
	where $ L(z)=L_1(z) +L_2(z)  $, $ \psi(s)=\psi_1(s)+ \psi_2(s) $ and $ K= k_1 +k_2 $.\vspace{1.5mm}
	
	\noindent{\bf Case 2.} Suppose that $ e^{h_2(z)-h_1(z)} $ is a constant. Clearly, $ h_2(z)-h_1(z) $ must be a constant. In fact, we see that $ h_2(z+c)-h_1(z+c) $ is also a constant. Thus it follows from \eqref{eq-3.5} that $ e^{h_1(z+c)-h_1(z)} $ is a constant, otherwise, we will arrive at a contradiction. Therefore, we conclude that $ h_1(z) = L_{21}(z) + k_3 $ and $ h_2(z) = L_{21}(z) + k_4 $, where $ L_{21}(z)= \gamma_1 z_1+ \gamma_2 z_2 +\cdots+\gamma_n z_n $ and $ \gamma_1, \ldots, \gamma_n, k_3, k_4\in\mathbb{C} $. Thus, from \eqref{eq-3.3} it follows that
	\begin{align*}
		g(z)=2(\gamma_1 z_1+ \gamma_2 z_2 +\cdots+\gamma_n z_n) + B, \;\;\mbox{where}\;\; B= k_3 + k_4.
	\end{align*}
	However, we see that the solution $f$ takes the following form
	\begin{align*}
		f(z)=\dfrac{Q_1(z)e^{L_{21}(z) + k_3}+Q_2(z)e^{L_{21}(z) + k_4}} {2\sqrt{a}} =\beta(z)e^{L_{21}(z)},
	\end{align*}
	where $ \beta(z) $ satisfies 
	\begin{align*}
		a\beta^2(z)+bP^2(z)\left(a_1\beta(z+c)+a_0\beta(z)\right)^2=e^{B}Q(z).
	\end{align*} 
	This completes the proof.
\end{proof}
\subsection{Proof of Theorem \ref{th-2.2}}
\begin{proof}[\bf Proof of Theorem \ref{th-2.2}] The equation \eqref{Meq-1.2} can be written as 
	\begin{align*}
		\left(\frac{\sqrt{a}f(z+c)}{e^{\frac{g(z)}{2}}}\right)^2+\left(\frac{\sqrt{b}P(z)\left(\frac{\partial f(z)}{\partial z_i}\right)}{e^{\frac{g(z)}{2}}}\right)^2=Q(z).
	\end{align*}
	By the similar argument being used in case of proof of Theorem \ref{th-2.1}, it can be easily shown that
	\begin{align}\label{eq-3.10}
		\begin{cases}
			f(z+c)=\dfrac{Q_1(z)e^{h_1(z)}+Q_2(z)e^{h_2(z)}}{2\sqrt{a}}\vspace{2mm}\\
			\dfrac{\partial f(z)}{\partial z_i}=\dfrac{Q_1(z)e^{h_1(z)}-Q_2(z)e^{h_2(z)}}{2i\sqrt{b}P(z)},
		\end{cases}
	\end{align}
	where
	\begin{align}\label{eq-3.11}
		h_1(z)=\dfrac{g(z)}{2}+\alpha(z) \;\;\;\;\;\mbox{and}\;\;\;\;\; h_2(z)=\dfrac{g(z)}{2}-\alpha(z).
	\end{align}
	Combining both the equations in \eqref{eq-3.9}, a simple computation shows that
	\begin{align}\label{eq-3.12}
		&\sqrt{a}Q_1(z+c) e^{h_1(z+c)-h_1(z)} -\sqrt{a}Q_2(z+c)e^{h_2(z+c)-h_1(z)} \\&\nonumber -i\sqrt{b}P(z+c)\left(Q_2(z)\frac{\partial h_2}{\partial z_i}+\frac{\partial Q_2}{\partial z_i}\right)e^{h_2(z)-h_1(z)}= i\sqrt{b}P(z+c)\left(Q_1(z)\frac{\partial h_1(z)}{\partial z_i} +\frac{ \partial Q_1}{\partial z_i}\right).
	\end{align}
	However, equation \eqref{eq-3.12} can be expressed as 
	\begin{align}\label{eq-3.13}
		g_1(z) + g_2(z) + g_3(z) = 1,
	\end{align} where
\begin{align}\label{Eq-4.18}
	\begin{cases}
		g_1(z)=\dfrac{\sqrt{a}Q_1(z+c)}{i\sqrt{b}P(z+c)\left(Q_1(z)\frac{\partial h(z)}{\partial z_i}+\frac{\partial Q_1}{\partial z_i}\right)}e^{h_1(z+c)-h_1(z)}\vspace{2mm}\\
		g_2(z)=-\dfrac{\sqrt{a}Q_2(z+c)}{i\sqrt{b}P(z+c)\left(Q_1(z)\frac{\partial h(z)}{\partial z_i}+\frac{\partial Q_1}{\partial z_i}\right)}e^{h_2(z+c)-h_1(z)}\vspace{2mm}\\
		g_3(z)=-\dfrac{Q_2(z)\frac{\partial h_2}{\partial z_i}+\frac{\partial Q_2}{\partial z_i}}{Q_1(z)\frac{\partial h_2}{\partial z_i}+\frac{\partial Q_1}{\partial z_i}}e^{h_2(z)-h_1(z)}.
	\end{cases}
\end{align} 
To complete the proof, it suffices to discuss the following  cases:\vspace{1.5mm}
	
\noindent{\bf Case I.} If $ e^{h_2(z)-h_1(z)} $ is non-constant, then it follows that the polynomial $ h_2(z)-h_1(z) $ must be non-constant. We claim that $ g_2 $ and $ g_3 $ are non-constants. Otherwise, if $g_2(z)=C_2$ and $g_3(z)=C_3$, where $ C_2 $ and $ C_3 $ are constants, then from \eqref{Eq-4.18}, it follows that
\begin{align*}
\begin{cases}
-\dfrac{\sqrt{a}Q_2(z+c)}{i C_2\sqrt{b}P(z+c)\left(Q_1(z) \frac{\partial h(z)}{\partial z_i}+\frac{\partial Q_1}{\partial z_i}\right)}=e^{h_1(z)-h_2(z+c)}\vspace{2mm}\\ \mbox{and} \;\;\;\;
-\dfrac{Q_2(z)\frac{\partial h_2}{\partial z_i}+\frac{\partial Q_2}{\partial z_i}}{C_3\left(Q_1(z)\frac{\partial h_2}{\partial z_i}+\frac{\partial Q_1}{\partial z_i}\right)}=e^{h_1(z)-h_2(z)}.
\end{cases}
\end{align*}
 It is easy to see that right side is a transcendental entire function but left side is a rational functions; both these equations leads us  to a contradiction.\vspace{1.2mm} 
 
 In view of Lemma \ref{lem-3.3}, we must have $ g_1(z)\equiv 1 $ \textit{i.e.,} we have 
	\begin{align}\label{eq-3.14}
		\dfrac{\sqrt{a}Q_1(z+c)}{i\sqrt{b}P(z+c)\left(Q_1(z)\frac{\partial h(z)}{\partial z_i}+\frac{\partial Q_1}{\partial z_i}\right)}\equiv e^{h_1(z)-h_1(z+c)}.
	\end{align}
	Evidently, it follows from \eqref{eq-3.14} that $ h_1(z)-h_1(z+c) $ must be constant and hence $ e^{h_1(z)-h_1(z+c)} $ becomes a constant, say $ A_{11} $, otherwise, we arrive at a contradiction. Assume that 
	\begin{align}\label{eq-3.15}
		\dfrac{\sqrt{a}Q_1(z+c)}{i\sqrt{b}P(z+c)\left(Q_1(z)\frac{\partial h(z)}{\partial z_i}+\frac{\partial Q_1}{\partial z_i}\right)} = A_{11}.
	\end{align} 
	In view of \eqref{eq-3.13} and \eqref{eq-3.14}, a simple computation shows that
	\begin{align}\label{eq-3.16}
		-\dfrac{\sqrt{a}Q_2(z+c)}{i\sqrt{b}P(z+c)\left(Q_2(z)\frac{\partial h_2}{\partial z_i}+\frac{\partial Q_2}{\partial z_i}\right)} =e^{h_2(z) -h_2(z+c)}.
	\end{align}
	Similarly, it follows from \eqref{eq-3.16} that $ h_2(z)-h_2(z+c) $ must be a constant, Hence, $ e^{h_2(z)-h_2(z+c)} $ must be a constant, say $ A_{12} $. Suppose that 
	\begin{align}\label{eq-3.17}
		-\dfrac{\sqrt{a}Q_2(z+c)}{i\sqrt{b}P(z+c)\left(Q_2(z)\frac{\partial h_2}{\partial z_i}+\frac{\partial Q_2}{\partial z_i}\right)} = A_{12}.
	\end{align}
	Rewriting \eqref{eq-3.15}, we get
	\begin{align}\label{eq-3.18}
		P(z+c)=\frac{\sqrt{a}Q_1(z+c)}{A_{11}i\sqrt{b}\left(Q_1(z)\frac{\partial h_1}{\partial z_i}+\frac{\partial Q_1}{\partial z_i}\right)}
	\end{align}
	which shows that $ Q_1(z) $ and $ \frac{\partial h_1}{\partial z_i} $ must be  constant, say $ \alpha_1 $ and $ \beta_1 $, respectively. Then $ h_1(z)= \beta_1 z_i+\phi(w) $, where $ w=(z_1, \ldots, z_{i-1}, z_{i+1}, \ldots, z_n) $. Thus, it follows from \eqref{eq-3.17} that $ P(z+c)=\frac{\sqrt{a}\alpha_1} {i\sqrt{b}A_{11}\alpha_1\beta_1}=P(z) $. Moreover, from \eqref{eq-3.17}, we obtain
	\begin{align*}
		P(z+c)=-\frac{\sqrt{a}Q_2(z+c)}{i\sqrt{b}A_{12}\left(Q_2(z)\frac{\partial h_2}{\partial z_i}+\frac{\partial Q_2}{\partial z_i}\right)}
	\end{align*}
	and from this, we conclude that $ Q_2(z) $ and $ \frac{\partial h_2}{\partial z_i} $ are non-zero constant, say $ \alpha_2 $ and $ \beta_2 $, respectively. Thus $ h_2(z)=\beta_2 z_i + \psi(w) $ and this shows that $ P(z+c)=-\frac{\sqrt{a}\alpha_2}{i\sqrt{b}\alpha_2 \beta_2 A_{12}}=P(z). $ Consequently, we have
	\begin{align*}
		\frac{\sqrt{a}\alpha_1}{i\sqrt{b} \alpha_1\beta_1 A_{11}}=-\frac{\sqrt{a}\alpha_2}{i\sqrt{b}\alpha_2\beta_2 A_{12}}
	\end{align*} 
	and this implies that $ A_{11}\beta_1+A_{12}\beta_2=0 $. Moreover, $ e^{h_1(z)-h_1(z+c)} =A_{11} $ implies that $ h_1(z)-h_1(z+c) =\ln A_{11} +2k\pi i $; $ k\in\mathbb{Z} $. That is 
	\begin{align*}
		\phi(w+c)-\phi(w)=-(\beta_1 c_i + \ln A_{11} +2k\pi i) = \mbox{constant}.
	\end{align*} 
	Hence, $ \phi(w)= a_1 z_1 +\cdots + a_{i-1} z_{i-1} + a_{i+1} z_{i+1} +\cdots + a_n z_n + r_1 $. Therefore, $ h_1(z) = a_1 z_1 +\cdots + a_{i-1} z_{i-1} + \beta_1 z_i + a_{i+1} z_{i+1} +\cdots + a_n z_n + r_1 $. Similarly, we obtain that $ h_2(z)= b_1 z_1 +\cdots + b_{i-1} z_{i-1} + \beta_2 z_i + b_{i+1} z_{i+1} +\cdots + b_n z_n + r_2 $. Consequently, 
	\begin{align*}
		\beta_2 e^{(a_1 -b_1) c_1 +\cdots + (a_{i-1} -b_{i-1}) c_{i-1} + (\beta_1 -\beta_2) c_i + (a_{i+1} -b_{i+1}) c_{i+1} +\cdots + (a_n -b_n) c_n} + \beta_1 = 0.
	\end{align*}
	\noindent{\bf Case-II:} If $ e^{h_2(z)-h_1(z)} $ is a constant, then $ h_2(z)-h_1(z) $ must be a constant. Then by \eqref{eq-3.13}, $ e^{h_1(z+c)-h_1(z)} $ is also a constant, otherwise, we get contradiction. Hence, we deduce that $ h_1(z)=L_1(z) + H(s)+ r_3 $ and $ h_2(z)=L_1(z) + H(s) + r_4 $, where $L_1(z)=\xi_1 z_1+ \cdots + \xi_n z_n$, $H(s)$ is a polynomial in $s=d_1 z_1 +\cdots+d_n z_n$ such that $d_1 c_1 +\cdots+d_n c_n=0$ with $H(z+c)=H(z)$ and $ r_3, r_4 \in\mathbb{C} $. Therefore, from \eqref{eq-3.11}, it is easy to see that
	\begin{align*}
		g(z) = 2L_1(z) + 2H(s) + r_5, \;\mbox{where}\;
		r_3 + r_4 =r_5.
	\end{align*} 
	Indeed, we have the form of the solution as
	\begin{align*}
		f(z)&=\frac{Q_1(z-c)e^{L_1(z-c)+ H(s) + r_3} + Q_2(z-c)e^{L_1(z-c)+ H(s)+ r_4}}{2\sqrt{a}} \\&=\gamma(z-c)e^{L_1(z-c)+ H(s)},
	\end{align*}
	where $ \gamma(z) $ satisfies 
	\begin{align*}
		a\gamma^2(z) + bP^2(z) e^{-2L_1(c)} \left(\frac{\partial \gamma(z-c)}{\partial z_i} + (\xi_i +d_i) \gamma(z-c)\right)^2 = e^{r_5}Q(z).
	\end{align*} 
	This completes the proof.
\end{proof}
\subsection{Proof of Theorem \ref{th-2.3}}
\begin{proof}[\bf Proof of Theorem \ref{th-2.3}]
Suppose that $ f(z) $ is a transcendental entire solution of finite order of the equation \eqref{Meq-1.3}. First of all, we write \eqref{Meq-1.3} as
\begin{align}\label{eq-3.19}
	aF^2 + 2\omega F G + bG^2 =1,
\end{align}
where $ F $ and $ G $ are defined by 
\begin{align}\label{eq-3.20}
	F = \dfrac{f(z)}{e^{\frac{g(z)}{2}}} \;\;\;\; \mbox{and} \;\;\;\; G= \dfrac{\gamma_1f(z + c) +\gamma_2 f(z)} {e^{\frac{g(z)}{2}}}
\end{align}
Further, we can transform \eqref{eq-3.1} to be 
\begin{align*}
	(\sqrt{a} F -\omega_1\sqrt{b} G)(\sqrt{a} F -\omega_2\sqrt{b} G) = 1.
\end{align*}
Since $ f $ is a finite order transcendental entire function and $ g $ is a polynomial, by Lemmas \ref{lem-3.1} and \ref{lem-3.2}, we see that there exists a polynomial $ p $ in $ \mathbb{C}^n $ such that
\begin{align}\label{eq-3.21}
	\sqrt{a} F -\omega_1\sqrt{b} G = e^{p} \;\;\;\;\mbox{and} \;\;\;\; \sqrt{a} F -\omega_2\sqrt{b} G =e^{-p}.
\end{align}
Using \eqref{eq-3.20} and \eqref{eq-3.21}, a simple computation shows that
\begin{align}\label{eq-3.22}
	f(z) = \dfrac{\omega_2 e^{p(z)} - \omega_1 e^{-p(z)}}{\sqrt{a}(\omega_2 -\omega_1)} e^{\frac{g(z)}{2}}, 
\end{align}
\begin{align}\label{eq-3.23}
	\gamma_1f(z + c) +\gamma_2 f(z) = \dfrac{e^{p(z)} - e^{-p(z)}}{\sqrt{b}(\omega_2 -\omega_1)} e^{\frac{g(z)}{2}}.
\end{align}
For brevity, let 
\begin{align}\label{eq-3.24}
	h_1(z)= \dfrac{g(z)}{2} + p(z) \;\;\;\;\mbox{and}\;\;\;\; h_2(z)= \dfrac{g(z)}{2} - p(z).
\end{align} 
The equations \eqref{eq-3.22} and \eqref{eq-3.23} can be written as
\begin{align}\label{eq-3.25}
	f(z) = \dfrac{\omega_2 e^{h_1(z)} - \omega_1 e^{h_2(z)}}{\sqrt{a}(\omega_2 -\omega_1)}
\end{align}
\begin{align}\label{eq-3.26}
	\gamma_1f(z + c) +\gamma_2 f(z) = \dfrac{e^{h_1(z)} - e^{h_2(z)}}{\sqrt{b}(\omega_2 -\omega_1)}.
\end{align}
Combining \eqref{eq-3.25} and \eqref{eq-3.26}, an easy computation shows that
\begin{align}\label{eq-3.27}
	R_{11} e^{h_1(z) -h_1(z+c)} + R_{12} e^{h_2(z) -h_1(z+c)} + R_{13} e^{h_2(z+c) -h_1(z+c)} \equiv 1,
\end{align}
where
\begin{align*}
	R_{11}=-\dfrac{(\omega_2\gamma_2\sqrt{b}-\sqrt{a})}{\omega_2\gamma_1 \sqrt{b}}, \;\; R_{12}=\dfrac{(\omega_1\gamma_2 \sqrt{b} -\sqrt{a})}{\omega_2\gamma_1 \sqrt{b}}\;\;\mbox{and}\;\; R_{13}=\dfrac{\omega_1}{\omega_2}.
\end{align*}
\noindent{\bf Case A:} If $ e^{h_2(z+c) -h_1(z+c)} $ is a constant, then $ h_2(z+c) -h_1(z+c) $ must be a constant. We set $ h_2(z+c) -h_1(z+c)=k $, where $ k\in\mathbb{C} $. It follows from \eqref{eq-3.24}, that $ p(z) $ is a constant. Let $ \xi=e^{p(z)} $, then from \eqref{eq-3.22} and \eqref{eq-3.23}, it is easy to see that
\begin{align}\label{eq-3.28}
	f(z) = B_{1} e^{\frac{g(z)}{2}} \;\; \mbox{and}\;\; \gamma_1f(z + c) +\gamma_2 f(z)= B_{2} e^{\frac{g(z)}{2}},
\end{align}
where 
\begin{align*}
B_{1} = \dfrac{\omega_2 \xi - \omega_1 \xi^{-1}}{\sqrt{a}(\omega_2 -\omega_1)} \;\;\mbox{and}\;\; B_{2}= \dfrac{\xi - \xi^{-1}}{\sqrt{b}(\omega_2 -\omega_1)}.
\end{align*}
It is easy see that $ B_1\neq 0 $ and $ B_2\neq 0 $. From \eqref{eq-3.28}, we see that 
\begin{align}\label{eq-3.29}
\frac{\gamma_1 B_1}{(B_2 -\gamma_2 B_1)}e^{\frac{g(z+c)-g(z)}{2}} =1.
\end{align}
Since $ g(z) $ is a polynomial, then \eqref{eq-3.29} implies $ g(z +c) - g(z)= \eta $, where $ \eta $ is a constant in $ \mathbb{C} $. Therefore, it follows that $ g(z)= L(z) + H(s) + B_3 $, where $ L(z)= a_1 z_1 +\cdots + a_n z_n $ and $ H(s) $ is a polynomial in $ s:= d_1 z_1 +\cdots + d_n z_n $ in $ \mathbb{C}^n $ with $ d_1 c_1 +\cdots + d_n c_n = 0 $ with $ H(z+c)= H(z) $, $ B_3\in\mathbb{C} $. Substituting $ g(z) $ into \eqref{eq-3.29} a simple computation shows that
\begin{align*}
	e^{\frac{a_1 c_1 +\cdots + a_n c_n}{2}} =\frac{(\omega_1 \gamma_2 \sqrt{b} -\sqrt{a}) -(\omega_2\gamma_2\sqrt{b} -\sqrt{a})\xi^2} {\gamma_1\sqrt{b}(\omega_2\xi^2 -\omega_1)}.
\end{align*} 
Therefore, in view of first equation of \eqref{eq-3.28}, we obtain
\begin{align*}
	f(z) = \frac{\omega_2\xi^2 -\omega_1}{\xi\sqrt{a}(\omega_2 -\omega_1)} e^{\frac{1}{2}(L(z) + H(s) + B_3)}.
\end{align*}
	
\noindent{\bf Case B:} If $ e^{h_2(z+c) -h_1(z+c)} $ is not a constant, then obviously, $ R_{11}\equiv 0 $ and $ R_{12}\equiv 0 $ cannot hold at the same time. Otherwise, from \eqref{eq-3.27} we see that $ R_{13} e^{h_2(z+c) -h_1(z+c)}\equiv 1 $, a contradiction. If $ R_{11}\equiv 0 $ and $ R_{12}\not\equiv 0 $, then in view of \eqref{eq-3.27}, we obtain
\begin{align}\label{eq-3.30}
	R_{12} e^{h_2(z) -h_1(z+c)} + R_{13} e^{h_2(z+c) -h_1(z+c)} \equiv 1.
\end{align}
Since $ e^{h_2(z+c) -h_1(z+c)} $ is not a constant, it follows that $ e^{h_2(z) -h_1(z+c)} $ is not a constant. Moreover, $ e^{h_2(z+c) -h_2(z)} $ is not a constant. Otherwise, $ h_2(z+c) -h_2(z)=\xi_1 $, where $ \xi_1\in\mathbb{C} $. From \eqref{eq-3.30}, we see that $ (R_{12} e^{-\xi_1} + R_{13}) e^{h_2(z+c) -h_1(z+c)} \equiv 1 $, which is a contradiction as $ e^{h_2(z+c) -h_1(z+c)} $ is non-constant. Hence, \eqref{eq-3.30} can be written as 
\begin{align}\label{eq-3.31}
	R_{12} e^{h_2(z)} + R_{13} e^{h_2(z+c)} - e^{h_1(z+c)} \equiv 0.
\end{align}
In view of Lemma \ref{lem-3.4}, from \eqref{eq-3.31}, we get a contradiction. Similarly, if $ R_{11}\not\equiv 0 $ and $ R_{12}\equiv 0 $, we get a contradiction. Thus, we conclude that $ R_{11}\not\equiv 0 $ and $ R_{12}\not\equiv 0 $.\vspace{1.2mm}

Since $ h_1(z), h_2(z) $ are polynomials and $ e^{h_2(z+c) -h_1(z+c)} $ is non-constant, then by Lemma \ref{lem-3.3} and with the help of \eqref{eq-3.27}, we obtain
\begin{align*}
	R_{11} e^{h_1(z) -h_1(z+c)}\equiv 1 \;\;\;\mbox{or}\;\;\;  R_{12} e^{h_2(z) -h_1(z+c)}\equiv 1.
\end{align*}
\noindent{\bf Sub-case B1:} Assume that $ R_{11} e^{h_1(z) -h_1(z+c)} \equiv 1 $. From \eqref{eq-3.27}, it yields that $ -\frac{R_{12}}{R_{13}} e^{h_2(z) -h_2(z+c)}\equiv 1 $. Since $ h_1(z), h_2(z) $ are polynomials, we see that $ h_1(z) -h_1(z+c)=\lambda_1 $ and $ h_2(z) -h_2(z+c)= \lambda_2 $, where $ \lambda_1, \lambda_2\in\mathbb{C} $. Thus, it follows that $ h_1(z)=L_1(z) + H_1(s) + D_1 $ and $ h_2(z)=L_2(z) +H_2(s) + D_2 $, where $ L_{l}(z)= a_{l1} z_1 +\cdots + a_{ln} z_n $ and $ H_{l}(s) $ for $ l=1, 2 $ are polynomial in $ s:= d_1 z_1 +\cdots + d_n z_n $ in $ \mathbb{C}^n $, $ d_1 c_1 +\cdots + d_n c_n = 0 $ with $ H_{l}(z+c)= H_{l}(z) $ for $ l=1, 2 $, and $ D_1, D_2\in\mathbb{C} $. Obviously, $ L_1(z) +H_1(s)\neq L_2(z) +H_2(s) $. Otherwise, we have $ h_2(z+c) -h_1(z+c) $ is a constant, which shows that $ e^{h_2(z+c) -h_1(z+c)} $ is a constant, a contradiction. Substituting $ h_1(z) $ and $ h_2(z) $ into $ R_{11} e^{h_1(z) -h_1(z+c)} \equiv 1 $ and $ -\frac{R_{12}}{R_{13}} e^{h_2(z) -h_2(z+c)}\equiv 1 $, we obtain
\begin{align*}
	-\frac{(\omega_2\gamma_2\sqrt{b}-\sqrt{a})}{\omega_2\gamma_1 \sqrt{b}} e^{-L_{1}(c)}\equiv 1 \;\;\;\;\mbox{and} \;\;\;\; -\frac{(\omega_1\gamma_2\sqrt{b}-\sqrt{a})}{\omega_1\gamma_1\sqrt{b}} e^{-L_{2}(c)}\equiv 1.
\end{align*}
Therefore, from \eqref{eq-3.7} it follows that $f$ takes the form
\begin{align*}
	f(z) = \dfrac{\omega_2 e^{L_1(z) +H_1(s) + D_1} - \omega_1 e^{L_2(z) +H_2(s) + D_2}}{\sqrt{a}(\omega_2 -\omega_1)}.
\end{align*}
Moreover, from \eqref{eq-3.6} we see that
\begin{align*}
	g(z) = h_1(z) + h_2(z) = L(z) + H(s) + D,
\end{align*}
where $ L(z):= L_1(z) + L_2(z) $, $ H(s):= H_1(s) + H_2(s) $ and $ D:= D_1 + D_2 $. \vspace{1.2mm}
	
\noindent{\bf Sub-case B2:} Assume that $ R_{12} e^{h_2(z) -h_1(z+c)} \equiv 1 $. Then, from \eqref{eq-3.27} it is easy to see that $ -\frac{R_{11}}{R_{13}} e^{h_1(z) -h_2(z+c)}\equiv 1 $. Since $ h_1(z), h_2(z) $ is a polynomial, it follows that $ h_2(z) -h_1(z+c)=\lambda_3 $ and $ h_1(z) -h_2(z+c)=\lambda_4 $, where $ \lambda_3, \lambda_4 \in \mathbb{C} $. A simple computation shows that $ h_1(z+2c) -h_1(z)= -\lambda_3 -\lambda_4 $ and $ h_2(z+2c) -h_2(z)= -\lambda_3 -\lambda_4 $. Thus, we deduce that $ h_1(z)= L(z) + H(s) + D_3 $ and $ h_2(z)= L(z) + H(s) + D_4 $, where $ L(z)= a_1 z_1 +\cdots + a_n z_n $ and $ H(s) $ is a polynomial in $ s:= d_1 z_1 +\cdots + d_n z_n $ in $ \mathbb{C}^n $ with $ d_1 c_1 +\cdots + d_n c_n = 0 $ with $ H(z+c)= H(z) $, and $ D_3, D_4\in \mathbb{C} $. Now, we see that $ h_2(z+c) -h_1(z+c)= D_4 - D_3 $, which shows that $ e^{h_2(z+c) -h_1(z+c)} $ is a constant, a contradiction. This completes the proof.
\end{proof}
\section{Declaration}
\vspace{1.6mm}

\noindent\textbf{Compliance of Ethical Standards:}\\

\noindent\textbf{Conflict of interest.} The authors declare that there is no conflict  of interest regarding the publication of this paper.\vspace{1.5mm}

\noindent\textbf{Data availability statement.}  Data sharing is not applicable to this article as no datasets were generated or analyzed during the current study.

\end{document}